\documentclass[11pt,psfig]{article}
\usepackage{amssymb,amsmath,amscd}
\usepackage{tikz-cd}
\usepackage[normalem]{ulem}
\usepackage[all]{xy}
\newtheorem{conj}{Conjecture}
\newtheorem{remark}{Remark}[section]
\newtheorem{lem}{Lemma}[section]
\newtheorem{thm}{Theorem}[section]
\newtheorem{cor}{Corollary}[section]
\newtheorem{prop}{Proposition}[section]
\newtheorem{defn}{Definition}

\newcommand{\f}[1]{\mathfrak{#1}}

\newcommand{\SPAN}{{\rm span}}

\newcommand{\mb}{\mathbb}
\newcommand{\ca}{{\#}}
\newcommand{\commentout}[1]{}

\newcommand{\mc}{\mathcal}
\newcommand{\arr}[1]{\left( \begin{array}{clcr} #1 \end{array} \right)}

\newcommand{\diag}{{\,\rm diag}}
\newcommand{\Hom}{{\rm {Hom}\,}}
\newcommand{\sgn}{\, {\rm sgn}}

\newcommand{\im}{{\rm Im}}
\newcommand{\cl}{{\rm cl}}
\newcommand{\supp}{{\rm supp}}

\begin{document}
\title{On the Gan-Gross-Prasad Conjecture for $U(p,q)$}
\author{Hongyu He  \\
Department of Mathematics \\
Louisiana State University \\
email: hongyu@math.lsu.edu \\
}
\date{}
\maketitle

\abstract{ In this paper, we give a proof of the Gan-Gross-Prasad conjecture for the discrete series of $U(p,q)$. There are three themes in this paper: branching laws of a small $A_{\f q}(\lambda)$, branching laws of discrete series and inductive construction of discrete series. These themes are linked together by a reciprocity law and the notion of invariant tensor product.}

\section*{Introduction}
In \cite{gp}, Gross and Prasad formulated a number of conjectures regarding the restrictions  of generic representations of the special orthogonal groups over a local field. These conjectures related the restriction problem to local root numbers.  In \cite{ggp}, Gan, Gross and Prasad extended these local conjectures to all classical groups. These local conjectures are known as the local Gan-Gross-Prasad (GGP) conjectures, or Gross-Prasad conjectures. Among the GGP conjectures, there was a very specific and interesting conjecture about the branching law of the discrete series representations for the real groups. Recently, there has been rapid development concerning the local GGP conjectures over the non-Archimedean fields starting with the work of Waldspurger (\cite{w}). As we understand, all cases of the non-Archimedean local conjectures are  close to being completely proved, with some standard assumptions.  For the Archimedean fields,  Gross and Wallach gave a proof of the Gross-Prasad conjecture for a class of small discrete series representation of $SO(2n+1)$.  Since then there has not been much progress towards the Gan-Gross-Prasad conjecture over the real numbers. The purpose of this paper is to give a proof of the GGP conjecture for the discrete series representations of $U(p,q)$. We shall also mention the recent work of Zhang that dealt with the global Gan-Gross-Prasad conjecture for the unitary group (\cite{zh}).\\
\\
Discrete series of $U(p,q)$ are parametrized by Harish-Chandra parameters. Following \cite{gp}, let $(\chi, z)$ be a Harish-Chandra parameter for $U(p,q)$ where $\chi \in \mathbb R^{p+q}$ is a sequence of distinct integers or half integers and $z \in \{\pm 1\}^{p+q}$ is a sequence of $+$ and $-$ corresponding to each entry in $\chi$. Here the total number of $+$'s must be $p$ and the total number of $-$'s must be $q$.  One may also interpret $z$ as a $(p,q)$-partition of $\chi$. Let $D(\eta, t)$ be a discrete series representation of $U(p-1, q)$. The Gan-Gross-Prasad conjecture gave a precise description of those $D(\eta, t)$ that appear as subrepresentations of $D(\chi, z)|_{U(p-1,q)}$. In addition, the multiplicity of these $D(\eta, t)$ must be all one.  Since the discrete spectrum of $D(\chi, z)|_{U(p-1,q)}$ only involves the discrete series, GGP conjecture produces a complete description of the discrete spectrum of $D(\chi, z)|_{U(p-1,q)}$. The recent work of Sun and Zhu settled the multiplicity one part generally for all irreducible unitary representations of $U(p,q)$ (\cite{sz}). \\
\\
To be more precise, GGP conjecture predicts that $D(\eta, t)$ appears in $D(\chi, z)|_{U(p-1,q)}$ if and only if $(\eta, t)$ and $(\chi, z)$ interlace each other in a very specific way. To describe this interlacing relation, let us recall the branching law for the compact group $U(p)$. Let $V_{\lambda}$ be an irreducible representation of $U(p)$ with highest weight $\lambda$ and $V_{\mu}$ be an irreducible representation of $U(p-1)$ with highest weight $\mu$. Suppose that $\lambda$ and $\mu$ are both arranged in descending order. Then $V_{\mu}$ appears in the branching law of the restriction of $V_{\lambda}$ if and only if $\lambda$ and $\mu$ satisfy the Cauchy interlacing relation:
$$\lambda_1 \geq \mu_1 \geq \lambda_2 \geq \mu_2 \geq \ldots \geq \lambda_{p-1} \geq \mu_{p-1} \geq \lambda_p.$$
This interlacing relation, expressed in terms of {\it Harish-Chandra parameters}, becomes a strict interlacing relation:
$$\chi_1  > \eta_1 > \chi_2 > \eta_2 > \ldots > \chi_{p-1} > \eta_{p-1} > \chi_p.$$
The interlacing relation predicted by GGP conjecture is a natural generalization of the Cauchy interlacing relation.
\begin{defn} We say that two Harish-Chandra parameters $(\chi, z)$ and $(\eta, t)$, of $U(p,q)$ and $U(p-1,q)$ respectively, satisfy the Gan-Gross-Prasad or simply the Gross-Prasad interlacing relation,  if one can line up $\chi$ and $\eta$ in the descending ordering such that the corresponding sequence of signs from $z$ and $t$ only has the following eight adjacent pairs
$$(\oplus +), (+ \oplus), (- \ominus), (\ominus -), (+-), (-+), (\oplus \ominus), (\ominus \oplus).$$
Here $\oplus$ and $\ominus$ represent $+1$ and $-1$ in $t$, and $+$ and $-$ represent $+1$ and $-1$ in $z$.
We call such a sign sequence the (interlacing) sign pattern of $(\chi, z)$ and $(\eta, t)$.
\end{defn}
Clearly, when there is neither $-$ nor $\ominus$, this interlacing relation is exactly the classical  Cauchy interlacing relation. The local GGP conjecture can then be reformulated as follows:
\begin{conj}[ \cite{gp}, \cite{ggp}]  
$D(\eta, t)$ appears as a subrepresentation of $D(\chi, z)|_{U(p-1,q)}$ if and only if $(\chi, z)$ and $(\eta, t)$ satisfy the GGP interlacing relation. 
\end{conj}
  We shall make a few remarks here.
  \begin{enumerate}
  \item In the literature, except in  \cite{gw}, the GGP conjectures are stated in terms of equivariant homomorphisms of smooth representations (\cite{gp}, \cite{ggp}). For discrete series, our version is equivalent to the smooth version.
  \item Originally, Gross and Prasad  stated their conjecture (Conjecture 12.27 \cite{gp}) for the orthogonal groups in terms of noncompactness of certain root basis associated with the Harish-Chandra parameters, {\it not in terms of the interlacing relation}.  Our formulation turns out to be equivalent to theirs. We can see the equivalence of these two versions as follows. For unitary groups, following \cite{gp},  consider the group $U(p,q, Q_1) \times U(p-1,q, Q_2)$ embedded in $U(p+q, q+p-1, Q_1-Q_2)$  where $Q_1, Q_2$ are the quadratic forms defining $U(p,q)$ and $U(p-1,q)$. Then the Harish-Chandra parameter $(\chi \oplus \eta, z \oplus t)$ defines a (positive) root basis of $U(p+q, q+p-1)$. Conjecture 12.27 says that $D(\eta, t)$ appears as a subrepresentation of $D(\chi, z)|_{U(p-1,q)}$ if and only if this root basis consists of only {\it noncompact} roots in $U(p+q, q+p-1)$. Clearly,  the noncompact roots from $U(p,q)$ and $U(p-1,q)$
 correspond exactly to adjacent signs 
 $$ (+-), (-+), (\oplus \ominus), (\ominus \oplus),$$
 the other noncompact roots from $U(p+q, q+p-1)$
  correspond exactly to adjacent sign
  $$(\oplus +), (+ \oplus), (- \ominus), (\ominus -).$$

\item In the framework of Langlands classification, the discrete series are grouped together in L-packets. For the discrete series of $U(p,q)$, each L-packet contains exact $\frac{ (p+q) !}{p! q !}$ discrete series representations. What GGP conjecture predicts and implies  is that at most one representation in each $L$-packet of $U(p-1,q)$ can occur as a subrepresentation in a fixed $D(\chi, z)$. In addition, each $D(\eta, t)$ can only appear as a subrepresentation in at most one discrete series representation in each L-packet of $U(p,q)$. In some sense, discrete series representations in the same $L$-packet naturally repel each other! This is quite remarkable. This phenomena seems to persist for tempered L-packets, as predicted by the GGP multiplicity one conjecture.

\end{enumerate}

 The main result in this paper confirms the Gan-Gross-Prasad conjecture for the discrete series.
 \begin{thm} The discrete spectrum of $D(\chi, z)|_{U(p-1,q)}$:
$$D(\chi, z)|_{U(p-1,q)}^{dis} = \hat{\oplus}_{(\eta, t)} D(\eta, t)$$
where the direct sum is taking over all those Harish-Chandra parameters $(\eta, t)$ such that $(\chi, z)$ and $(\eta, t)$ satisfy the GGP interlacing relation.
 \end{thm}
Recently, Beuzart-Plessis seems to have proved the GGP multiplicity one conjecture for tempered L-packets (\cite{bp}). In our context, his result says that for a fixed $(\chi, \eta)$, there exists a unique $(z, t)$ such that $D(\eta, t) $ is a subrepresentation of $D(\chi, z)$ for a suitable pair of unitary groups. In comparison, our result not only implies the multiplicity one conjecture, but also pins down exactly the $(z, t)$ parameter, namely $(\chi, z)$ and $(\eta, t)$ must satisfy the GGP interlacing relation.   Beuzart-Plessis's method is quite different from ours, as we do not use the trace formula in our paper.  We now describe very briefly the main ideas of this paper. \\
 \\
Recall that discrete series representations, at least for classical groups, belong to a larger class of unitary representations, known as representations with nonzero cohomology. Representations with nonzero cohomology are very important in the theory of automorphic forms. They contribute to the cohomolgy of Shimura varieties. In \cite{vz}, Vogan and Zuckerman gave a characterization of representations with nonzero cohomology.  These representations can be constructed using Zuckerman's derived functor  (see for example, \cite{kv}), as $A_{\f q}(\lambda)$'s. Therefore, they are also known as $A_{\f q}(\lambda)$.   The branching law of $A_{\f q}(\lambda)$ becomes important, partly because of its implications in the theory of automorphic forms. This is, perhaps, one of the motivations of the GGP conjectures. In any case, we shall refer the reader to the recent survey article by Kobayashi for the branching laws of $A_{\f q}(\lambda)$(\cite{ko2}).\\
\\
 One focus of this paper is indeed the branching law of some small $A_{\f q}(\lambda)$. Let $U(n,n)$ be the isometry group of a Hermitian form with signature $(n,n)$ on $\mathbb C^{2n}$.  Let  us  consider $A_{\f q}(\lambda)$ of $U(n,n)$ where $\f q$ has the real Levi factor $U(r,s) \times U(s,r)$ with $r+s=n$ and $\lambda$ is weakly dominant with respect to $\f q$.  We denote them by $A_{r,s;s,r}(k_1, k_2)$ with $k_1 \geq k_2$ normalized parameter (see 2.3 \cite{pt}).  As we will show later in this paper, many branching laws of more general $A_{\f q}(\lambda)$ are extrinsically connected with the branching laws of $A_{r,s; s,r}(k_1, k_2)$. \\
 \\
Let us decompose $\mathbb C^{2n}$ into a direct sum of two orthogonal subspaces such that the Hermitian forms restricted to these two subspaces are all nondegenerate. Let $G_1 \times G_2$ be the subgroup of $U(n,n)$ that preserves this decomposition. We call such a group diagonally embedded in $U(n,n)$. Let $U(p,q) \times U(q,p)$ with $p+q=n$ be a subgroup diagonally embedded in $U(n,n)$. The first theorem we proved is about the branching law of $A_{r,s;s,r}(k_1, k_2)|_{U(p,q) \times U(q,p)}$. We state the case for $p+q$ even, when $k_1, k_2$ are  integers. For $p+q$ odd, see Prop \ref{discrete:equalcase2}.
\begin{thm}\label{equal}  Let $p+q=r+s=n$ be an even positive integer. Let $k_1 \geq k_2$ be integers. Then the discrete spectrum of $A_{r,s;s,r}(k_1, k_2)|_{ U(p,q) \times U(q, p)}$ is the direct sum of
$$[D(\eta, -t^{\prime})^* \otimes {\det}^{k_2}] \boxtimes [D(\eta, t) \otimes {\det}^{k_1}]$$
where
\begin{eqnarray}
& \#\{ \eta_i >0 \mid t_i =1 \}+ \#\{  \eta_i < 0 \mid t_i =-1 \}   =  r ,\\
& \#\{\eta_i >0 \mid  t_i =-1 \}+ \#\{  \eta_i < 0  \mid t_i =1 \}  =   s, \\
& \#\{ \eta_i \in ( k_2-k_1, 0) \mid t_i=1 \}  =  \# \{   \eta_i \in ( k_2-k_1, 0) \mid t_i=-1 \} , \\
& t^{\prime}_i  =  \left\{ \begin{array}{cc} 
					t_i & \mbox{ if $\eta_i \notin ( k_2-k_1,0)$} \\
                  -t_i & \mbox{ if  $\eta_i \in (k_2-k_1, 0)$}.
                  \end{array}
                  \right.
\end{eqnarray} 
Here $\pi^*$ stands for the contragredient representation of $\pi$, $\eta_i$ will all be half integers.   
\end{thm}
This theorem is proved by connecting these $A_{\f q}( \lambda)$ with Howe's correspondence (\cite{howe}) through a theorem due to A. Paul and P. Trapa (\cite{pt}).  The main technique is applying various properties of invariant tensor products defined in Section 1. The description of Howe's correspondence in the equal rank case, due to Li and Paul, is crucial here (\cite{li1}, \cite{p1}). \\
\\
Let $U(p-1, q) \times U(q+1, p)$ be another diagonally embedded subgroup of $U(n,n)$. Using similar ideas, we establish a second branching law:
$$A_{r,s;s,r}(k_1,k_2)|^{dis}_{U(p-1,q) \times U(q+1, p)} \supseteq \oplus D(\chi^{\prime}, z^{\prime}) \boxtimes D(\chi, z),$$
Here $(\chi, z)$ satisfies a set of equations similar to Theorem. \ref{equal}, and $(\chi^{\prime}, z^{\prime})$ corresponds to $(\chi, z)$ in a definite way. Some limit of discrete series of $U(q+1,p)$ will appear in the discrete spectrum. Hence we use $\supseteq$ to indicate this fact. 
  See section 5 for the details. \\
\\
With the discrete spectrum of these two branching laws in hand, it is not hard to prove  that every discrete series representation will appear in some $A_{r,s;s,r}(k_1, k_2)$. See Thm. \ref{discrete:equalcase} and Cor. \ref{induction1}.  We need the following lemma of reciprocity to link these two branching laws together.
\begin{lem}[Reciprocity] Let $H_1$ and $H_2$ be two separable locally compact Hausdorff topological group, and let $\pi$ be a unitary representation of $H_1 \times H_2$. Then for all irreducible unitary representations $\sigma_1$ of $H_1$ and $\sigma_2$ of $H_2$, there are canonical isometric isomorphisms
$${\Hom}_{H_1}(\sigma_1, {\Hom}_{H_2}(\sigma_2, \pi)) \cong {\Hom}_{H_1 \times H_2}(\sigma_1 \hat{\otimes} \sigma_2, \pi) \cong {\Hom_{H_2}}(\sigma_2, {\Hom}_{H_1}(\sigma_1, \pi))$$
of Hilbert spaces. Here $\Hom_{H_i}$ refers to continuous $H_i$-equivariant operators for Hilbert spaces.
\end{lem}
We would like to thank the second referee for suggesting this form of reciprocity.\\
\\
Now consider 
$$G= U(p+q, p+q); \qquad \pi=A_{r,s;s,r}(k_1, k_2); $$
$$\begin{array}[c]{ccc}
H_2= U(q,p) &\stackrel{}{\subseteq}& G_2= U(q+1,p)\\
\updownarrow\scriptstyle{}&&\updownarrow\scriptstyle{}\\
G_1 =U(p,q) &\stackrel{}{\supseteq} & H_1=U(p-1,q) 
\end{array}$$
where $G_1 \times H_2$ and $G_2 \times H_1$ are diagonally embedded in $U(p+q, p+q)$. By the two branching laws and reciprocity, we will have
\begin{equation}
\begin{split}
 & \Hom_{U(p-1,q)}(D(\chi^{\prime}, z^{\prime}), D(\eta^{\prime}, -t^{\prime})^* \otimes {\det}^{k_2}) \\
 = &\Hom_{U(p-1,q)}(D(\chi^{\prime}, z^{\prime}), \Hom_{U(q,p)}(D(\eta, t)\otimes {\det}^{k_1}, \pi))\\
  =& \Hom_{U(q,p)}(D(\eta, t)\otimes {\det}^{k_1}, \Hom_{U(p-1,q)}(D(\chi^{\prime}, z^{\prime}), \pi))\\
  = & \Hom_{U(q,p)}(D(\eta, t)\otimes {\det}^{k_1}, D(\chi, z)) \\
\end{split}
\end{equation}
This will allow us to establish branching laws of discrete series inductively. Now  every discrete series will appear as a subrepresentation of some $A_{r,s;s,r}(k_1,k_2)$. Hence all discrete series representations will be covered in the induction process. In addition, only discrete series representations can appear in the discrete spectrum of $D(\chi, z)|_{U(q,p)}$ (see Theorem \ref{l2}). Hence we will be able to give a complete description of 
$D(\chi, z)|^{dis}_{U(q,p)}$ inductively. The rest of the proof involves some detailed numerical and combinatorial analysis on sequences and sign patterns. Surprisingly, the multiplicity one theorem comes for free from the induction process. \\
\\
\commentout{
We should remark that our approach can be extended to give a complete branching law, both continuous spectrum and discrete spectrum. We must prove a theorem analogous to Theorem \ref{} for continuous spectrum. Then the tempered version of GGP conjecture is expected to follow when the infinitesimal characters are nonsingular.}
We shall point out a remarkable feature of our proof. Unlike the proof of GGP conjectures for the non-Archimedean local field,  our proof given here is extrinsic by nature, namely, the discrete series representation of $U(p,q)$ is studied while being embedded as a subrepresentation of $A_{r,s;s,r}(k_1, k_2)$. In other words, we use the symmetries \lq\lq outside\rq\rq  the group $U(p,q)$ to establish our result. In the absence of an intrinsic proof, we may ask,  to what extent, we can understand a representation through extrinsic methods. Clearly, there is  Howe's theory of dual reductive pair which can be regarded as an extrinsic studies of representations (\cite{howe}).  The question is then, whether there are  representations other than the Weil representation crucial in the extrinsic studies of representations. Undoubtedly, the  $A_{\f q}(\lambda)$ in the \lq\lq middle dimension \rq\rq will deserve serious consideration. Roughly these are the $A_{\f q}(\lambda)$'s associated with various complex Siegel parabolic subalgebra $\f q$ or its analogies. For $GL(2n, \mathbb R)$, these $A_{\f q}(\lambda)$ are known as the Speh representations. In the case of $U(n,n)$, all of them can be obtained by Howe's correspondence up to a central character. We take full advantage of this connection in our analysis of these $A_{\f q}(\lambda)$'s. However, for other groups, one cannot obtain all these middle dimensional $A_{\f q}(\lambda)$ through theta lifting of one dimensional representations.  Hence, different ideas will be needed to study the branching law of the $A_{\f q}(\lambda)$ in middle dimension.  For symplectic groups or orthogonal groups, we expect the GGP conjectures can be proved once these branching laws are established. We make some conjectures in this context in Section 8.\\
\\ 
Finally, we would like to thank Prof. Gross for encouraging us to work on the Gan-Gross-Prasad conjecture and thank Wee-Teck Gan for updating us the current status in the real cases and $p$-adic cases. We also want to thank the referees for their very detailed comments and suggestions.
\subsection{Conventions}
\uline{All the Hilbert spaces and topological groups in this paper are assumed to be separable}. \uline{All groups are assumed to be locally compact, Hausdorff, topological group}. Unless otherwise stated, the statements we underline in this paper will be assumed through out the rest of this paper. 
We call a number half integer if and only if it is an integer plus $\frac{1}{2}$. The constant $C$ is used as a positive symbolic constant. 

\section{Invariant Tensor Products}
 In this section, we shall formulate the concept of invariant tensor product. Our motivation of defining invariant tensor product comes from developing the analytic theory of Howe's correspondence (\cite{li}\cite{theta}) and constructing unipotent representations (\cite{hen} \cite{hearthur}). In this paper, we shall study the invariant tensor product associated with smooth representations equipped with  invariant Hermitian forms.  This is slightly different from the view point of \cite{hen} and \cite{hearthur}. \\
 \\
Let $V$ be a linear vector space over $\mathbb C$. Define the complex vector space $V^c$ by letting $V =V^c$ as a real vector space and the complex  multiplication $ \lambda v \in V^c$ to be $\overline{\lambda} v$ in $V$.  Let $G$  be a  group. Let $(\pi, V)$ be a linear representation of $G$ over $\mathbb C$. \uline{ Here a linear representation refers to a representation of an abstract group, without any consideration of the topology}. Define the representation $(\pi^c, V^c)$ by  $\pi^c(g)=\pi(g)$. If $G$ is a topological group and $(\pi, V)$ is a unitary representation, then $(\pi^c, V^c) \cong (\pi^*, V^*)$ by Riesz representation theorem. Here $(\pi^*, V^*)$ is the (unitary) contragredient representation of $(\pi, V)$.

\begin{defn}[Invariant Tensor Product]\label{itp} Let $G$ be a locally compact topological group and $d g$ be a left invariant Haar measure. Let $(\pi, H_{\pi})$ and $(\pi_1, H_{\pi_1})$ be two unitary representations of  $G$. Let $V$ and $V_1$ be two dense subspaces of $H_{\pi}$ and $H_{\pi_1}$. Assume that $\forall \ \  u,v \in V, u_1, v_1 \in V_1$, 
 $$
 \int_{G} (\pi(g)v, u)(\pi_1(g) v_1, u_1) d g$$
 converges absolutely. Define the averaging operator
$$\mathcal A: V \otimes V_1 \rightarrow \Hom (V^c \otimes V_1^c, \mathbb C)$$
by
 \begin{eqnarray}~\label{average}
\mathcal A( v \otimes v_1)(u \otimes u_1) & = & \int_{G} {\big(} (\pi \otimes \pi_1)(g) (v \otimes v_1), (u \otimes u_1){\big)}  d g \\
& = & \int_{G} (\pi(g)v, u)(\pi_1(g) v_1, u_1) d g 
\end{eqnarray}
and extending this by linearality to $V \otimes V_1$. Here $\Hom$ refers to algebraic homomorphisms. 
Define the invariant tensor product $V \otimes_{G} V_1$ to be the image   $\mathcal A(V \otimes V_1)$.  Whenever we use the notation $V \otimes_{G} V_1$, {\bf we assume $V \otimes_{G} V_1$ is well-defined}, that is, the integral (~\ref{average}) converges absolutely for all $u,v \in V, u_1, v_1 \in V_1$. Denote $\mathcal A(v \otimes v_1)$ by $v \otimes_{G} v_1$. 
\end{defn}
\begin{remark}  For $G$ compact, various forms of the space $V \otimes_G V_1$ have long been used implicitly in the literature. For example, when $G$ is compact and $V$ irreducible ( necessarily finite dimensional), the space $\mc H_{\pi} \otimes_G \mc H_{\pi_1}$ is always well-defined. We have $\mc H_{\pi} \otimes_G \mc H_{\pi_1}= \Hom_G(\mc H_{\pi}^c, \mc H_{\pi_1})$. Hence $\mc H_{\pi} \otimes_G \mc H_{\pi_1}$ can be used to compute the multiplicity of $\mc H_{\pi}^c$ in $\mc H_{\pi_1}$. When $G$ is noncompact, then $V \otimes_G V_1$ is still related to a certain multiplicity space in a very delicate manner. For $G$ a classical group,  $V \otimes_G V_1$ appeared implicitly in \cite{li}, in the theory of doubling zeta integrals of Piatetski-Shapiro and Rallis, and in other texts.
\end{remark}
\begin{remark}
It is clear from our definition that we can always assume that $V$ and $V_1$ are $G$-stable without loss of any generalities. Indeed, let ${}^G V$ be the linear vector space spanned by $\pi(g) v (v \in V)$. If $V \otimes_G V_1$ is well-defined, then ${}^G V \otimes_G {}^G V_1$ will also be well-defined and vice versa.
\end{remark}

\begin{defn}[Hermitian form]\label{hf} In addition, suppose that $G$ is unimodular. Then
$$(v \otimes_G v_1)(u \otimes u_1)=(u \otimes_G u_1)(v \otimes v_1)= \int_{G} (\pi(g)v, u)(\pi_1(g) v_1, u_1) d g.$$
Define
$$(v \otimes_G v_1, u \otimes_G u_1)_G= \int_{G} (\pi(g)v, u)(\pi_1(g) v_1, u_1) d g.$$
This form defines a non-degenerate pairing between $V \otimes_G V_1$ and $V^c \otimes_G V_1^c$.
It yields a Hermitian form on $V \otimes_G V_1$.
\end{defn}
\subsection{Basic Properties and Equivalence of Representations}
Given two linear representations $V$ and $W$ of $G$, if $V$ and $W$ are equivalent as linear representations of $G$,  we write $V \cong W$ or $V \cong_G W$. 
If $G$ is trivial, $V \cong W$ simply means that $V$ is isomorphic to $W$. 
Invariant tensor products have the following properties.
\begin{lem}\label{basicproperty} Let $\mc H$, $\mc H_1$ and $\mc H_2$ be unitary representations of $G$. Let $V$, $V_1$ and $V_2$ be dense subspaces of $\mc H$, $\mc H_1$ and $\mc H_2$ respectively. Then
\begin{enumerate}
\item (commutativity) $V \otimes_G V_1 \cong V_1 \otimes_G V.$
\item (associativity I) $ V \otimes_G (V_1 \otimes V_2) \cong (V \otimes V_1) \otimes_G V_2$.
\item (associativity II) If $G$ acts trivially on $V_2$, then $ V \otimes_G (V_1 \otimes V_2) \cong (V \otimes_G V_1) \otimes V_2$.
\item Regarding $(\, , \,)_G$ as a Hermitian form on $V \otimes V_1$, let $\mc R$ be the radical of this form. Then $\ker(\mc A)=\mc R$ and $V \otimes_G V_1 \cong V \otimes V_1 / \mc R$.
\item Suppose that $V$ and $V_1$ are  $G$-invariant subspaces of $\mc H_{\pi}$ and $\mathcal H_{\pi_1}$ respectively. Then  $V^c$ and $V_1^c$ will be  $G$-invariant subspaces of $\mc H_{\pi}^c$ and $\mathcal H_{\pi_1}^c$ respectively.  We have
$$V \otimes_G V_1 \subseteq \Hom_{G}( V^c \otimes V_1^c, \mathbb C).$$
\end{enumerate}
\end{lem} 
The equivalences here are all canonical linear isomorphisms without consideration of topology. Defining topology is a little bit subtle. We shall confine ourselves to the case where the Hermitian structure  is positive definite. Then the topology can be defined in terms of the Hermitian structure. \\
\\
Let $(\pi_1, V_1)$ and $(\pi_2,V_2)$ be two continuous representations of $G$ on two pre-Hilbert spaces such that the group actions preserve the pre-Hilbert structures respectively.
Then $(\pi_i, V_i)$ completes to a unitary representation $(\pi_i, \mc H_i)$. \uline{ We say that $(\pi_1, V_1)$ is equivalent to $(\pi_2, V_2)$ if $(\pi_1, \mc H_1)$ is equivalent to $(\pi_2, \mc H_2)$. By abusing notation, we denote this by $V_1 \cong V_2$ or $\pi_1 \cong \pi_2$. To specify the group $G$, we may write $V_1 \cong_G V_2$ or $\pi_1 \cong_G \pi_2$. This notation applies to all continuous representations equipped with $G$-invariant inner products}. 
\begin{defn} Let $V$ and $V_1$ be two $G$-representations equipped with nondegenerate $G$-invariant inner product such that the group actions are continuous.  Definitions \ref{itp} and \ref{hf} apply.
\end{defn}
Throughout this paper, we mostly work in the category of smooth representations. 
It is more convenient to apply this definition. Of course when we complete $V$ and $V_1$  to unitary representations, this new definition coincides with Definitions \ref{itp} and \ref{hf}.

\subsection{Constructing Representations}

\begin{defn}\label{rep} Let $G_1$ and $G_2$ be two  groups. Let $(\pi, \mc H_{\pi})$ be a unitary representation of 
$G_1 \times G_2$ and $(\pi_1, \mc H_{\pi_1})$ be a unitary representation of $G_1$. Let $V$ be a dense $G_2$-invariant subspace of $\mc H_{\pi}$. Let $V_1$ be a dense subspace of $\mc H_{\pi_1}$. If 
$V \otimes_{G_1} V_1$ is well-defined, we define:
$$(\pi \otimes_{G_1} \pi_1) ( g_2) (u \otimes_{G_1} u_1)= \pi(g_2) u \otimes_{G_1} u_1 \qquad (X_2 \in \f g_2, u \in V, u_1 \in V_1).$$
Then $(\pi \otimes_{G_1} \pi_1, V \otimes_{G_1} V_1)$ is a linear representation of $G_2$.
\end{defn} 
We shall remark that if $V$ in Lemma \ref{basicproperty} has a $H$ action that commutes with the $G$ action, then all equivalences in Lemma \ref{basicproperty} hold as equivalences of linear $H$-representations.\\
\\
Notice that $\pi(g_2)$ acts on $V^c \otimes V_1^c$ via $V^c$. It induces a contragredient action on $\Hom(V^c \otimes V_1^c, \mathbb C)$. We can  also view the action of $(\pi \otimes_{G_1} \pi_1 )(g_2)$ as this contragredient action  restricted to the subspace $V \otimes_{G_1} V_1$. 
\begin{lem}~\label{reprn}  Suppose that $G_1$ is unimodular. Then the Hermitian form $(\, , \,)_{G_1}$ on $V \otimes_{G_1} V_1$ is $G_2$-invariant.
\end{lem}
Proof: Let $u, v \in V; u_1, v_1 \in V_1$ and $g_2 \in G_2$. Write $\sigma= \pi \otimes_{G_1} \pi_1$. Then
\begin{equation}
\begin{split}
 & (\sigma(g_2) (u \otimes_{G_1} u_1), v \otimes_{G_1} v_1)_{G_1}  \\
 = & \int_{G_1} (\pi(g_1) \pi(g_2) u, v)(\pi_1(g_1) u_1, v_1) d g_1 \\
 = & \int_{G_1} (\pi(g_1) u, \pi(g_2^{-1}) v) (\pi_1(g_1) u_1, v_1) d g_1 \\
 = & (u \otimes_{G_1} u_1, \pi(g_2^{-1}) v \otimes_{G_1} v_1)_{G_1} \\
 = & (u \otimes_{G_1} u_1, \sigma(g_2^{-1})( v \otimes_{G_1} v_1))_{G_1}. 
 \end{split}
\end{equation}
So $(\, , \,)_{G_1}$ is $G_2$-invariant. 
We have thus shown that $(,)_{G_1}$ on $ V \otimes_{G_1}  V_1$ is an invariant Hermitian form with respect to $(\pi \otimes_{G_1} \pi_1) (G_2)$. $\Box$\\
\\
When the Hermitian form $(\, , \,)_{G_1}$ on $V \otimes_{G_1} V_1$ is positive definite, we equip $V \otimes_{G_1} V_1$ with the pre-Hilbert space structure defined by $(\, , \,)_{G_1}$.
\subsection{Associativity}
The main theorem proved in this section is as follows.
\begin{thm}[Associativity] \label{asso}
In the setting of Definition \ref{rep}, suppose that $G_1$ and $G_2$ are both unimodular. Let $(\pi_2, \mc H_2)$ be a unitary representation of $G_2$ with a dense subspace $V_2$ such that 
\begin{enumerate}
\item $V_2 \otimes_{G_2} V$ is well-defined.
\item the functions
$$(\pi_2(g_2)u_2, v_2)(\pi(g_1 g_2) u, v)(\pi_1(g_1) u_1, v_1) \qquad (u_2, v_2 \in V_2; u, v \in V; u_1, v_1 \in V_1)$$
are all in $L^1(G_1 G_2)$.
\end{enumerate}  
Suppose that for $i=1,2$, the Hermitian form on $V_i \otimes_{G_i} V$ is positive definite and the linear representation $\pi_i \otimes_{G_i} \pi$ completes to a unitary representation. Then
$$(V_2 \otimes_{G_2} V) \otimes_{G_1} V_1 \cong V_2 \otimes_{G_2} ( V \otimes_{G_1} V_1).$$
\end{thm}
Proof: First of all, we have
$$V_2 \otimes_{G_2} (V \otimes_{G_1} V_1) \subseteq \Hom (V_2^c \otimes (V \otimes_{G_1} V_1)^c, \mathbb C) \cong \Hom(V_2^c \otimes (V^c \otimes_{G_1} V_1^c), \mathbb C) $$
$$(V_2 \otimes_{G_2} V) \otimes_{G_1} V_1 \subseteq \Hom( (V_2 \otimes_{G_2} V)^c \otimes V_1^c, \mathbb C) \cong \Hom( (V_2^c \otimes_{G_2} V^c) \otimes V_1^c, \mathbb C).$$
Recall the averaging operators
$$\mc A_2: V_2^c \otimes V^c \rightarrow V_2^c \otimes_{G_2} V^c, \qquad \mc A_1: V^c \otimes V_1^c \rightarrow V^c \otimes_{G_1} V_1^c.$$
are surjective.
 We identify $\phi \in V_2 \otimes_{G_2} (V \otimes_{G_1} V_1)$ with the element $\tilde \phi$  in 
$\Hom( V_2^c \otimes V^c \otimes V_1^c, \mathbb C)$ as follows:
\commentout{\[
\begin{array}{ccc}
 V_2^c \otimes V^c \otimes V_1^c  & \xrightarrow{\mathcal A_1} & V_2^c \otimes (V \otimes_{G_1} V_1)^c & 
&\searrow{\tilde{\phi}}&\Big\downarrow\rlap{$\scriptstyle\phi$}& 
&& \mathbb C &.
\end{array}
\]}
\[
\begin{tikzcd}
 V_2^c \otimes V^c \otimes V_1^c \arrow{r}{\mathcal A_1} \arrow[swap]{dr}{\tilde \phi} & V_2^c \otimes (V \otimes_{G_1} V_1)^c 
 \arrow{d}{\phi} \\
& \mathbb C.
\end{tikzcd}
\]
Similarly we identify $\psi \in ( V_2 \otimes_{G_2} V )\otimes_{G_1} V_1$ with the element  $\tilde \psi$ in 
$\Hom( V_2^c \otimes V^c \otimes V_1^c, \mathbb C)$. By Fubini's theorem, for any $v_i, u_i \in V_i$ and
$v, u \in V$, we have
$$\widetilde{[v_2 \otimes_{G_2} (v \otimes_{G_1} v_1)]}(u_2, u, u_1)= \widetilde{[(v_2 \otimes_{G_2} v) \otimes_{G_1} v_1]}(u_2, u, u_1)$$
$$= \int_{G_1 \times G_2} (\pi_2(g_2) v_2, u_2) (\pi(g_1 g_2) v, u)(\pi_1(g_1)v_1, u_1) d g_1 d g_2. $$
We see that $(V_2 \otimes_{G_2} V) \otimes_{G_1} V_1$ and $ V_2 \otimes_{G_2} ( V \otimes_{G_1} V_1)$ can be identified as the same subspace of  $\Hom( V_2^c \otimes V^c \otimes V_1^c, \mathbb C)$. The associativity follows immediately. $\Box$

\section{Discrete Spectrum and Square Integrable Representations }

\underline{Let $G$ be a locally compact topological group with a left invariant Haar measure $d g$}.   
Let $\hat G$ be the space of equivalence classes of irreducible unitary representations of $G$ equipped with the Fell topology. \underline{From now on suppose that $G$ is CCR or equivalently liminaire (\cite{wallach} \cite{di})}.  Then for any unitary representation $(\pi, \mc H_{\pi})$, there exists a Borel measure $d_{\pi}$ on $\hat G$, and a Hilbert space $M_{\pi}(\sigma)$ for each $\sigma$ defined almost everywhere with respect to $d_{\pi}$, such that
$$\mc H_{\pi} \cong \int_{\sigma \in \hat G} \mc H_{\sigma} \hat{\otimes} M_{\pi}(\sigma)d_{\pi}(\sigma).$$
Here $G$ acts trivially on $M_{\pi}(\sigma)$.
If $d_{\pi}(\sigma) > 0$, then $\sigma$ appears as a subrepresentation of $\pi$. We say that $\sigma$ is in the discrete spectrum of $\pi$. We may simply write $\sigma \in \pi$. In this situation, $M_{\pi}(\sigma)$ is well-defined. We write
$\mc H_{\pi}(\sigma)$ to be the closure of  the direct sum of all subrepresentations of $\pi$ equivalent to $\sigma$. Then $\mc H_{\pi}(\sigma) \cong_G \mc H_{\sigma} \hat{\otimes} M_{\pi}(\sigma)$. We call $\mc H_{\pi}(\sigma)$ the $\sigma$-isotypic subspace of $\mc H_{\pi}$. We call the Hilbert space $M_{\pi}(\sigma)$ the multiplicity space. The dimension of $M_{\pi}(\sigma)$ is called the multiplicity of $\sigma$ in $\pi$, often denoted by $m_{\pi}(\sigma)$ or $m(\pi, \sigma)$. Notice that $\sigma$-isotypic subspace $\mc H_{\pi}(\sigma)$ is only well-defined for $\sigma$ in the discrete spectrum of $\pi$. We define $(\pi^{dis}, \mc H_{\pi}^{dis})$ to be the closure of the direct sum of all
$\mc H_{\pi}(\sigma)$. \\
\\
\underline{Suppose from now on that $G$ is unimodular}. Then $L^2(G)$ has a left regular $G$-action and a right regular $G$-action. These two actions commute. So $L^2(G)$ is a unitary  representation of $G \times G$. There is the abstract Plancherel theorem proved by I. Segal.
\begin{thm}[Segal]
Let $G$ be a locally compact unimodular CCR topological group. Then
there exists a Borel measure $d \mu$ on $\hat G$ such that $L^2(G)$ decomposes  as 
$$\int_{\sigma \in \hat G} \mc H_{\sigma} \hat{\otimes} \mc H_{\sigma^*} d \mu(\sigma)$$
and for any $f \in L^2(G) \cap L^1(G)$, we have
$$(f, f)=\int_{\hat G} Tr(\sigma(f)^* \sigma(f)) d \mu(\sigma).$$
\end{thm}
$d \mu$ is call the Plancherel measure. The discrete spectrum of $L^2(G)$ are often called discrete series representations.

\subsection{Square Integral Representations}
We call a unitary representation $(\pi, \mc H_{\pi})$ {\it square integrable} or $L^2$ if there is a dense subspace $V $ in $\mc H_{\pi}$ such that all matrix coefficients 
$$g \rightarrow (\pi(g) v_1, v_2) \qquad (v_1, v_2 \in V)$$
are in $L^2(G)$. Here we do not assume $\pi$ is irreducible. However if $\pi$ is irreducible, then all matrix coefficients will be in $L^2(G)$. The (equivalence classes of) irreducible square integrable representations are precisely those appearing in the discrete spectrum of $L^2(G)$ (\cite{di}). 
\begin{thm}\label{l2} Let $(\pi, \mc H_{\pi})$ be a square integrable representation of $G$. If $\sigma \in \hat G$ appears in the discrete spectrum of $\pi$, then $\sigma$ is square integrable.
\end{thm}
Therefore only discrete series can appear in the discrete spectrum of a square integrable representation. This result is known to expert, at least for real reductive groups. I have not been able to find any proofs or references in the literature. Because of the importance of this theorem in this paper, I will supply a proof here. The  main ingredients come from \cite{li1}. 
\begin{lem}[Mackey's Schur Lemma]
Let $G$ be a locally compact topological group. Let $(\pi, \mc H_{\pi})$ and $(\tau, \mc H_{\tau})$ be two unitary representations of $G$. Let $T$ be a closed $G$-equivariant ( unbounded) operator defined
on a dense subspace $V$ of $\mc H_{\pi}$, i.e., $T: V \rightarrow \mc H_{\tau}$.  Then there is a $G$-equivariant isometry between $(\ker T)^{\perp}$ and $\cl(\im(T))$.
\end{lem}
Here $\ker T$ is closed and the image of $T$ may not be closed; $\cl(\im(T))$ is the closure of the image of $T$ in $\mc H_{\tau}$. This version of Mackey's Schur lemma is quoted from \cite{fd}.\\
\\
Proof of Theorem \ref{l2}: Suppose that the linear subspace $V$ is dense in $\mc H_{\pi}$ such that the matrix coefficients for $V$ are all in $L^2(G)$. Without loss of generality, suppose that $V$ is $G$-stable. Fix $v \in V$. Consider
$$j_{v}: V \rightarrow L^2(G)$$
defined by
$$j_{v}(u)(g)=(u, \pi(g) v), \qquad (u \in V).$$
It is easy to verify that $j_{v}(\pi(h) u)(g) =j_{v}( u)(h^{-1} g)$ for any $h \in G$. So
$j_{v}$ intertwines the action of $\pi$ and the left regular action. Furthermore $j_{v}$ is closable (see Page 717 \cite{li1}). Let $j^0_{v}$  be the closure of $j_{v}$. If $u \in {\rm Dom}(j_{v}^0)$ then there exists a sequence $\{ u_i \} \subseteq {\rm Dom}(j_{v})$ such that $u_i \rightarrow u$ in $\mc H_{\pi}$ and $j_v(u_i)(g) \rightarrow j_v^0(u)(g)$ in $L^2(G)$.\\
\\
Let $u \in \ker(j_v^0)$. Then $j_v^0(u)(g)=(u, \pi(g) v) \equiv 0$. Hence
$$\ker ({j_v}^0) \subseteq \{ u \in \mc H_{\pi} \mid (u, \pi(g) v)=0 \, \, \forall \,  g \}.$$
It follows that 
$$\ker({j_v}^0)^{\perp} \supseteq \SPAN \{\pi(g) v : g \in G \}.$$
Therefore $\ker({j_v}^0)^{\perp}$ contains  (the closure of) the cyclic space generated by $v$. 
By Mackey's Schur lemma, $\ker({j_v}^0)^{\perp} \cong \cl(\im({j_v}^0))$. Since $\cl(\im(j_v^0)) \subseteq L^2(G)$,   $\sigma \in {\ker(j_v^0)^{\perp}}$ implies that $\sigma \in L^2(G)$.
We conclude that if $\sigma \in {\cl(\SPAN(\pi(g) v \mid g \in G))}$, then $\sigma \in L^2(G)$. \\
\\
Now suppose that $\sigma \in \pi$. Then there exists a $v \in V$ such that the orthogonal projection of $\mc H_{\pi}(\sigma)$ onto the $\cl(\SPAN(\pi(g) v \mid g \in G))$ is not zero.
Otherwise, $\mc H_{\pi}(\sigma)$ will be orthogonal to $V$, hence zero. In any case, we will have
$\sigma \in \cl(\SPAN(\pi(g) v \mid g \in G))$ for some $v$. Hence $\sigma$ is a discrete series representation.
$\Box$\\
\\
There are often continuous spectrum in $\pi$. Our definition of $\sigma \in \pi$ does not cover the continuous spectrum. We may use the notion that $\sigma \in_{wk} \pi$ if $\sigma$, as an irreducible unitary representation, is weakly contained in $\pi$. Then both continuous and discrete spectrum are counted for. In this paper, we will not discuss $\sigma \in_{wk} \pi$. If we were to discuss the continuous spectrum of the restriction of the discrete series, then we will need to use the notion of weak containment.
\subsection{Multiplicity Spaces}
Let $\pi$ be a square integrable representation of $G$. Let $\sigma \in \pi$. Then $\sigma$ is a discrete series representation. We are interested in the multiplicity space $M_\pi(\sigma)$.
\begin{thm}\label{multi} Let $\pi$ be a square integrable representation of $G$ with respect a dense subspace $V$. Let $\sigma$ be an irreducible square integrable representation of $G$. Let $W$ be any $G$-invariant subspace of $\mc H_{\sigma^c}$. Then $V \otimes_{G} W$ is well-defined and the associated Hermitian form is positive definite. 
\begin{enumerate}
\item $M_{\pi}(\sigma) \cong  \cl(V \otimes_G  W)$,
\item $\mc H_{\pi}(\sigma)\neq 0$ if and only if $V \otimes_G  W \neq 0$.
\item $\mc H_{\pi}(\sigma) \cong_G \mc H_{\sigma} \hat{\otimes} \, \cl(V \otimes_G  W)$,
\item $\mc H_{\pi}^{dis} \cong_G \hat{\oplus}_{\sigma \,  L^2} \mc H_{\sigma} \hat{\otimes} \cl(V \otimes_G W)$.
\end{enumerate}
In particular, all statements hold if  $W=\mc H_{\sigma^c}$ or $W=\mc H_{\sigma^c}^{\infty}$.
\end{thm}
Proof: Recall that $\sigma^c$ is defined on $\mc H_{\sigma}^c$ with the same group action as $\sigma$. The matrix coefficients of $\sigma^c$ are exactly the conjugates of the matrix coefficients of $\sigma$. 
\begin{enumerate}
\item Suppose that  $M_{\pi}(\sigma) = 0$. Then $\mc H_{\pi}(\sigma) =0$. We must have
$V \otimes_G  W = 0$.  The contrapositive statement of this is proved as Lemma 2.2 (b) in \cite{li1} (Pg. 917).  
\commentout{
Otherwise, for some $v \in V$, the map $j_v: V \rightarrow L^2(G)$ given in the proof of Theorem \ref{l2} has the property that the closure of its image
$\cl(j_v(V)) \subseteq L^2(G)$ contains $\sigma$ as a subrepresentation. So by Mackey's Schur lemma and proof of Theorem \cite{l2}, the closure of the cyclic space generated by $v$ contains $\sigma$ as a subrepresentation. So $\mc H_{\pi}(\sigma) \neq 0$. This is a contradiction. \\
}
We have one direction of (2) and (1)(3)(4) all follow from this in the case $M_{\pi}(\sigma) = 0$. 
 
\item Suppose that $M_{\pi}(\sigma) \neq 0$. Let $\{e_1, e_2, \ldots \}$ be an orthonormal basis of
$M_{\pi}(\sigma)$. Let $P_{\sigma}$ be the canonical projection from $\mc H_{\pi}$ onto $\mc H_{\pi}(\sigma)$. Let $P_{\sigma}^i$ be the canonical map $\mc H_{\pi} \rightarrow \mc H_{\sigma}$ given by first projecting $\mc H_{\pi}$ onto $\mc H_{\sigma} \otimes \mathbb C e_i$, the $i$-th copy of $\mc H_{\sigma}$ and then identifying it with $\mc H_{\sigma}$. For any $u \in V$ and $x \in W$, define
$$j(u \otimes x)=\sum_{i} (P_{\sigma}^i (u), x) e_i.$$
This is a well-defined map from $V \otimes W$ into $M_{\pi}(\sigma)$.
Since $\{ e_i \}$ is an orthonormal basis, for any $v \in V$, $y \in W$,
$$(j(u \otimes x), j(v \otimes y))=( \sum_{i} (P_{\sigma}^i (u), x) e_i, \sum_{i} (P_{\sigma}^i (v), y) e_i)= \sum_{i} (P_{\sigma}^i (u), x) \overline{(P_{\sigma}^i (v), y)}.$$
On the other hand, we have

\begin{equation}
\begin{split}
(u \otimes_G x, v \otimes_G y) &= \int_G (\pi(g) u, v)(\sigma^c(g) x, y) d g \\
 &= \int_G (\pi(g) u, v)\overline{(\sigma(g) x, y)} d g \\
 &= \int_G (\pi(g) P_{\sigma}(u), P_{\sigma}(v)) \overline{(\sigma(g) x, y)} d g \\
 &= \int_G \sum_{i} (\sigma(g) P_{\sigma}^i(u), P_{\sigma}^i(v)) \overline{(\sigma(g) x, y)} d g \\
 &= \sum_{i} \int_G (\sigma(g) P_{\sigma}^i(u), P_{\sigma}^i(v)) \overline{(\sigma(g) x, y)} d g \\
 &= \frac{1}{d_{\sigma}} \sum_i (P_{\sigma}^i (u), x) \overline{(P_{\sigma}^i (v), y)}.
\end{split}
\end{equation}
The third equation follows from the fact that all matrix coefficients of $\mc H_{\sigma} \hat{\otimes} M_{\pi}(\sigma)$ are square integrable.
Changing order of $\int$ and $\sum$ is valid because
$\sum_{i} (\sigma(g) P_{\sigma}^i(u), P_{\sigma}^i(v))$ converges to $(\pi(g) P_{\sigma}(u), P_{\sigma}(v))$ in $L^2(G)$. Here $d_{\sigma}$ is the formal degree of $\sigma$.\\
\\
We see that 
$(j(u \otimes x), j(v \otimes y))=d_{\sigma} (u \otimes_G x, v \otimes_G y).$ 
By multilinear algebra, $j$ induces a map from $V \otimes_G W$ into
$M_{\pi}(\sigma)$ and the kernel of this map is $\{0\}$. It follows that the Hermitian form associated with $V \otimes_G W$ is positive definite.\\
\\
Now we want to show that $j(V \otimes W)$ is dense in $M_{\pi}(\sigma)$. By definition, $j(u \otimes x)=\sum_{i} (P_{\sigma}^i (u), x) e_i=
\sum_i (u, (P_{\sigma}^i)^* x) e_i.$ Here $(P_{\sigma}^i)^*$ is the adjoint of $P_{\sigma}^i$ and it maps vectors in $\mc H_{\sigma}^c$ to $\mc H_{\pi}^c$. Suppose that there is $\sum_i a_i e_i \in M_{\pi}(\sigma)$ perpendicular to all $j(V \otimes W)$. Then
$$0=(j(u \otimes x), \sum_{i} a_i e_i)=\sum_i (u, (P_{\sigma}^i)^* x) \overline{a_i} =(u, \sum_{i} a_i (P_{\sigma}^i)^* x)$$ for all $u \in V$ and $x \in W$. Since $V$  is dense in $\mc H_{\pi}$,  $\sum_{i} a_i (P_{\sigma}^i)^* x=0$. Now for a fixed nonzero $x \in W \subseteq \mc H_{\sigma^c}$, $\{ (P_{\sigma}^i)^* x \}$ is an orthogonal set in $\mc H_{\pi}(\sigma)$. So $a_i \equiv 0$. It follows that the image of $j(V \otimes  W)$ is dense in $M_{\pi}(\sigma)$. (1) is proved.
Furthermore $V \otimes_G W \neq \{0\}$.  The other direction of (2) is proved. (3) and (4) follow immediately. 
\end{enumerate}
$\Box$\\

\begin{cor} Under the same hypothesis as Theorem \ref{multi}, we have
$$V \otimes_G  W \cong V \otimes_G \mc H_{\sigma^c}^{\infty} \cong V \otimes_G \mc H_{\sigma^c}.$$
\end{cor}

\begin{remark} Following the convention in 1.1, the equivalence here means the closures of the invariant tensor products are equivalent:
$$cl(V \otimes_G  W) \cong cl(V \otimes_G \mc H_{\sigma^c}^{\infty}) \cong cl(V \otimes_G \mc H_{\sigma^c}).$$
Indeed, the closure of the invariant tensor product $V \otimes_G W$ is independent of the choices of the  dense subspace $W$ of the irreducible unitary representation $\sigma^c$. All the invariant tensor product we will be dealing with will have this property.  This validates the  use of $\cong$ to mean equivalences at the closure level throughout this paper.
\end{remark}
\section{Howe's Correspondence and $A_{\f q}(\lambda)$}
Let $(G, G^{\prime})$ be a real reductive dual pair in a real symplectic group $Sp$ (\cite{howe}). Let $\tilde{Sp}$ be the metaplectic covering of $Sp$ and $\{1, \epsilon \}$ the preimage of the identity.
Let $\tilde G$ and $\tilde G^{\prime}$ be the preimage of $G$ and $G^{\prime}$ respectively. Let
$\omega$ be the oscillator representation of $\tilde{Sp}$ with a lowest weight vector. Let $\mc R(\tilde G, \omega)$ and $\mc R(\tilde G^{\prime}, \omega)$ be the equivalence classes of irreducible admissible  representations of $\tilde G$ and $\tilde G^{\prime} $ occurring as quotients of $\omega^{\infty}$. Howe proved that there is a one-to-one correspondence between
$\mc R(\tilde G, \omega)$ and $\mc R(\tilde G^{\prime}, \omega)$ (\cite{howe}). This correspondence is called Howe's correspondence. Howe's correspondence  has been shown to possess nice properties related to parabolic induction (\cite{ku}), cohomological induction (\cite{li1}) and unitarity (\cite{li}, \cite{unit}). \\
\\
It is easy to see that all representations in $\mc R(\tilde G, \omega)$ have the property that
$\pi(\epsilon)= -1$. Representations with this property are called genuine representations, in the sense that they do not descend into representations of $G$.
\subsection{Howe's Correspondence for $U(p,q)$}
Howe's correspondence for $U(p,q)$ was studied in \cite{p1} \cite{p2}. In the equal rank case, A. Paul gave an explicit description of Howe's correspondence in terms of Langlands parameters. In particular, in this case, discrete series representations correspond to discrete series representations. Here we follow loosely the notation from \cite{pt}.\\
\\
Consider the dual pair $(U(p,q), U(r,s))$. It can be constructed as a subgroup of $Sp_{2(p+q)(r+s)}(\mathbb R)$ as follows. Let $V=M(p+q,r+s, \mathbb C)$ be the real vector space of $p+q$ by $r+s$ complex matrices. For $X, Y \in V$, define a real symplectic form
$$\langle X, Y \rangle= \Im \left(( Tr \arr{I_p & 0 \\ 0 & -I_q} X \arr{I_r & 0 \\ 0 & -I_s} \overline{Y^t} \right).$$
\underline{Let $U(p,q)$ and $U(r,s)$ be the groups preserving the Hermitian forms defined by}
$$\arr{I_p & 0 \\ 0 & -I_q} \qquad \arr{I_r & 0 \\ 0 & -I_s}$$
respectively. Let $U(p,q)$
act on $V$ from left and $U(r,s)$ act on $V$ from right. Then $(U(p,q), U(r,s))$ becomes a dual reductive pair in $Sp(V)$. Under the metaplectic covering of $Sp(V)$, the preimages are
$$\tilde U(p,q)=\{ (\lambda, g) \mid g \in U(p,q), \lambda \in \mathbb C, \lambda^2=\det(g)^{r-s} \},$$
$$\tilde U(r,s)=\{ (\lambda, g) \mid g \in U(r,s), \lambda \in \mathbb C, \lambda^2=\det(g)^{p-q} \}.$$
Following \cite{pt}, we denote these two groups by $U(p,q)^{\sqrt{\det^{r-s}}}$ and $U(r,s)^{\sqrt{\det^{p-q}}}$. In particular, when $r-s$ is odd, $U(p,q)^{\sqrt{\det^{r-s}}}$ is connected and 
will be a double covering of $U(p,q)$; when $r-s$ is even, $U(p,q)^{\sqrt{\det^{r-s}}}$ is disconnected and 
will be simply  $U(p,q) \times \{1, \epsilon \}$. Hence $U(p,q)^{\sqrt{\det^{r-s}}}$ only depends on the parity of $r-s$. Since we will study Howe's correspondence across different $U(r,s)$'s, we denote $U(p,q)^{\sqrt{\det^{r-s}}}$ by $U(p,q)^o$ if $r-s$ is odd and by $U(p,q)^e$ is $r-s$ is even. \\
\\
For $r-s$ even, genuine irreducible representations of $U(p,q)^{\sqrt{\det^{r-s}}}$ are of the form $\pi \otimes \sgn$ with $\pi$ an irreducible representations of $U(p,q)$ and $\sgn$ the sign representation of $\{1, \epsilon \} \cong \mathbb Z_2$. We parametrize the genuine irreducible representations of $U(p,q)^e$ by the irreducible representations of $U(p,q)$. For $r-s$ odd, a genuine representation of $U(p,q)^{\sqrt{\det^{r-s}}}$ can always be obtained from a representation of $U(p,q)$ by tensoring with $\det^{\frac{1}{2}}$:
$${\det}^{\frac{1}{2}}( \lambda, g)=\lambda (\det g)^{-\frac{r-s}{2}+\frac{1}{2}} \qquad (\,\, \forall \,\,\, (\lambda, g) \in U(p,q)^{\sqrt{\det^{r-s}}}) .$$ 
Nevertheless, we parametrize the genuine representations, for example the discrete series, exactly the same way as other connected reductive groups.  The infinitesimal character of genuine irreducible representation of $U(p,q)^o$ will differ from those of $U(p,q)$ by a shift of a half integer.\\
\\
Let us denote Howe's correspondence from $\mc R(\tilde U(p,q), \omega)$ to $\mc R(\tilde U(r,s), \omega)$ by
$\theta(p,q;r,s)$. It often suffices to use $\theta(r,s)(\pi)$ since the group $\tilde U(p,q)$ can often be read from the representation $\pi$ itself. This seems to be the convention adopted by several authors. However, since we will be dealing with several dual pairs at the same time, it is necessary to use $\theta(p,q;r,s)$. {\bf If $\pi \notin \mc R(\tilde{ U}(p,q), \omega)$, then we write $\theta(p,q;r,s)(\pi) =0$}. We will also use $\omega(p,q;r,s)$ to denote the oscillator representation $\omega$ when the dual pair $(U(p,q), U(r,s))$ is considered. 
\subsection{Discrete Series representations}
Let $T$ be the maximal torus diagonally embedded in $G=U(p,q)$ or $G=U(p,q)^o$.  Let $\mathfrak{g}_{\mathbb C}$ be the complex Lie algebra of $G$. Let $\f h= i \f t$. Identify $\f h$ with $\mathbb R^{p+q}$. Fix once for all a positive root system for $(\f g_{\mathbb C}, \f h)$ such that the dominant Weyl chamber is parametrized by
$$a_1 \geq a_2 \geq \ldots \geq \ldots \geq a_{p+q}.$$
Let $\rho$ be the half sum of the positive roots. Then 
$$\rho=(\frac{p+q-1}{2}, \frac{p+q-3}{2}, \ldots, -\frac{p+q-3}{2}, -\frac{p+q-1}{2}).$$
The (equivalence classes of) discrete series representations of $U(p,q)$ are parametrized by the following data
\begin{enumerate}
\item
 $\lambda=(\lambda_1 > \lambda_2 \ldots > \lambda_{p+q}) \in \rho+\mathbb Z^{p+q}$,
 \item a partition of $\lambda$ into $\lambda^+ \in \mathbb R^p$ and $\lambda^- \in \mathbb R^q$.
 \end{enumerate}
 Genuine discrete series representations of $U(p,q)^o$ can then be parametrized by
 \begin{enumerate}
\item
 $\lambda=(\lambda_1 > \lambda_2 \ldots > \lambda_{p+q}) \in \rho+\mathbb Z^{p+q}+\mathbf{\frac{1}{2}}$,
 \item a partition of $\lambda$ into $\lambda^+ \in \mathbb R^p$ and $\lambda^- \in \mathbb R^q$.
 \end{enumerate}
 Here $\mathbf c$ denotes the constant vector $(c,c,\ldots c)$.\\
 \\
  We call $(\lambda^+, \lambda^-)$ the Harish-Chandra parameter. For the time being, we allow $\lambda^+$ and $\lambda^-$ to be in any order. They become unique up to the action of the Weyl group of the maximal compact group of $G$. We can take $\lambda$ to be the infinitesimal character. The $L$-packets are parametrized by $\lambda$ and representations in the $L$-packet are parametrized by the partition. We denote the discrete series representation with Harish-Chandra parameter $(\lambda^+, \lambda^-)$ by $D(\lambda^+, \lambda^-)$. \\
\\
For any $\lambda^{\pm}$, let $\lambda_{+}^{\pm}$ be the positive portion and $\lambda_{-}^{\pm}$ be the negative portion. Then $\lambda^{\pm}=(\lambda_{+}^{\pm}, \lambda_{-}^{\pm})$ if $0$ does not appear in $\lambda^{\pm}$; $\lambda^{\pm}=(\lambda_{+}^{\pm}, 0, \lambda_{-}^{\pm})$ if $0$ appears in $\lambda^{\pm}$. Notice $0$ can only appear in  $\lambda$ once. It is easy to see that
$$D(\lambda^+, \lambda^-)^* \cong D(-\lambda^+, -\lambda^{-}), \qquad D(\lambda^+, \lambda^-) \otimes {\det}^k \cong  D(\lambda^+ + {\mathbf k}, \lambda^-+ \mathbf{k}).$$
\subsection{Discrete Series and Howe's Correspondence}
We shall now summarize some results of J.-S. Li in \cite{li1} and A. Paul in \cite{p1} \cite{p2} concerning discrete series representations of $U(p,q)$ under Howe's correspondence in the equal rank case.
\begin{thm}[Li, Paul]\label{lipaul} Suppose that $r+s =p+q$. Then Howe's correspondence maps discrete series to discrete series or zero. Let $\sigma$ be a discrete series representation of  $G=U(p,q)^{\sqrt{\det^{r-s}}}$. Then $\omega(p,q; r,s)^{\infty} \otimes_G \mc H_{\sigma}^{\infty}$ is well-defined and is infinitesimally equivalent to $\theta(p,q; r,s)(\sigma^c)$. 
Howe's correspondence for discrete series is given by
$$\theta(p,q; r,s) (D((\lambda^+_{+}, \lambda_{-}^+), (\lambda_{+}^{-}, \lambda_{-}^{-})))=D((\lambda^{+}_{+}, \lambda_{-}^{-}),(\lambda_{+}^{-}, \lambda_{-}^{+})),$$
where $\lambda$'s are all half integers and 
$$\ca(\lambda^+_{+})+\ca( \lambda_{-}^+)=p, \qquad \ca(\lambda_{+}^{-})+\ca( \lambda_{-}^{-})=q, $$
$$ \ca(\lambda^{+}_{+})+\ca( \lambda_{-}^{-})=r, \qquad \ca(\lambda_{+}^{-})+\ca( \lambda_{-}^{+})=s.$$
No other discrete series occur in $\mc R(U(p,q)^{\sqrt{\det^{r-s}}}, \omega)$ or $\mc R(U(r,s)^{\sqrt{\det^{p-q}}}, \omega)$. 
\end{thm}
The equations regarding the cardinality  are not necessary and they are dictated by the size of the maximal compact groups. We attach these equations here out of abundant caution. The correspondence was established by Li and the exhaustion was given by Paul.\\
\\
We remark that for $p+q$ even, a Harish-Chandra parameter for $G=U(p,q)^e$ consists of half-integers since $\rho$ consists of half integers; for $p+q$ odd, a Harish-Chandra parameter for $U(p,q)^o$ consists of half integers due to the covering.  
We also remark that the closure of $\omega(p,q; r,s)^{\infty} \otimes_G \mc H_{\sigma^c}^{\infty}$ is the multiplicity space of $\sigma$. Hence $\sigma \boxtimes \theta(p,q;r,s)(\sigma)$ is in the discrete spectrum of $\omega|_{\tilde U(p,q) \times \tilde U(r,s)}$. \\
\\
Consider the spectrum of $\omega|_{\tilde U(p,q) \times \tilde U(r,s)}$. Assume that $p+q \leq r+s$. Then  $\omega|_{\tilde U(p,q)}$ is  $L^2$. Li showed that if a discrete series representation $\pi$ of $\tilde U(p,q)$ appears in the discrete spectrum of $\omega|_{\tilde U(p,q)}$, then it will appear in Howe's correspondence. Equivalently, if a discrete series $\pi$ of $\tilde U(p,q)$ does not appear in Howe's correspondence, then $\pi$ does not appears in the discrete spectrum of $\omega|_{\tilde U(p,q)}$. By Theorems \ref{lipaul} and \ref{multi}, we have the following
\begin{cor}[Discrete Spectrum of $\omega(p,q; r, s)$, equal rank case]
Suppose $p+q=r+s$.
The discrete spectrum, 
$$\omega|_{\tilde U(p,q)}^{dis}=\omega|_{\tilde U(p,q) \times  \tilde U(r,s)}^{dis}= \omega|_{\tilde U(r,s)}^{dis} \cong \hat{\oplus} [\sigma \boxtimes \theta(p,q;r,s)(\sigma)]$$
where $\sigma$ and $\theta(p,q;r,s)(\sigma)$ are
 given by the  recipe  in Theorem \ref{lipaul}.
\end{cor}
\begin{thm}[Li, Paul]\label{lipaul2} Suppose that $r+s=p+q+1$. Then Howe's correspondence maps discrete series of the larger group $\tilde U(r,s)$ to discrete series of the smaller group $\tilde U(p,q)$ or zero. In this situation, the correspondence  is given by
$$\theta(p,q; r,s) (D((\lambda^+_{+}, \lambda_{-}^+), (\lambda_{+}^{-}, \lambda_{-}^{-})))=D((\lambda^{+}_{+}, 0, \lambda_{-}^{-}),(\lambda_{+}^{-}, \lambda_{-}^{+})),$$
where $\lambda$'s are all nonzero integers and 
$$\ca(\lambda^+_{+})+\ca( \lambda_{-}^+)=p, \qquad \ca(\lambda_{+}^{-})+\ca( \lambda_{-}^{-})=q, $$
$$ \ca(\lambda^{+}_{+})+\ca( \lambda_{-}^{-})+1=r, \qquad \ca(\lambda_{+}^{-})+\ca( \lambda_{-}^{+})=s;$$
or
$$\theta(p,q; r,s) (D((\lambda^+_{+}, \lambda_{-}^+), (\lambda_{+}^{-}, \lambda_{-}^{-})))=D((\lambda^{+}_{+},  \lambda_{-}^{-}),(\lambda_{+}^{-},0, \lambda_{-}^{+})),$$
where $\lambda$'s are all nonzero integers and 
$$\ca(\lambda^+_{+})+\ca( \lambda_{-}^+)=p, \qquad \ca(\lambda_{+}^{-})+\ca( \lambda_{-}^{-})=q, $$
$$ \ca(\lambda^{+}_{+})+\ca( \lambda_{-}^{-})=r, \qquad \ca(\lambda_{+}^{-})+\ca( \lambda_{-}^{+})+1=s.$$
The discrete series on the right hand side of the equations are all the discrete series occurring in $\mc R(U(r,s)^{\sqrt{\det^{p-q}}}, \omega)$.  Furthermore, when $\sigma$ is in the discrete series of $G=U(p,q)^{\sqrt{\det^{r-s}}}$, $\omega(p,q; r,s)^{\infty} \otimes_G \mc H_{\sigma}^{\infty}$ is well-defined and is infinitesimally equivalent to $\theta(p,q; r,s)(\sigma^c)$.
\end{thm}
 Generally speaking, in the case $r+s=p+q+1$, Howe's correspondence does not map the discrete series of the smaller group $\tilde{U}(p,q)$ exclusively to discrete series of the larger group $\tilde{U}(r,s)$. There are some discrete series representations  occurring in $\mc R(U(p,q)^{\sqrt{\det^{r-s}}}, \omega)$ and $\theta(p,q;r,s)$ maps these representations to  limits of discrete series of
 $U(r,s)^{\sqrt{\det^{p-q}}}$. The details are completely worked out by A. Paul.  Since the statement in Theorem \ref{lipaul2} is sufficient for the purpose of this paper, we will not give the complete statement of the correspondence for all discrete series in $\mc R(U(p,q)^{\sqrt{\det^{r-s}}}, \omega)$. The reader should notice that the discrete series in $\mc R(U(r,s)^{\sqrt{\det^{p-q}}}, \omega)$ must have a zero in their Harish-Chandra parameters. 
 
 \begin{cor}[Discrete Spectrum of $\omega(p,q; r, s)$:\,\,\, $r+s=p+q+1$ ]
Suppose $r+s=p+q+1$.
The discrete spectrum, 
$$\omega|_{\tilde U(p,q)}^{dis}=\omega|_{\tilde U(p,q) \times  \tilde U(r,s)}^{dis}= \omega|_{\tilde U(r,s)}^{dis}$$
contains both
$$ D((\lambda^+_{+}, \lambda_{-}^+), (\lambda_{+}^{-}, \lambda_{-}^{-})) \boxtimes D((\lambda^{+}_{+}, 0, \lambda_{-}^{-}),(\lambda_{+}^{-}, \lambda_{-}^{+}))$$
and 
$$ D((\lambda^+_{+}, \lambda_{-}^+), (\lambda_{+}^{-}, \lambda_{-}^{-})) \boxtimes D((\lambda^{+}_{+},  \lambda_{-}^{-}),(\lambda_{+}^{-},0, \lambda_{-}^{+})),$$ given by the  recipes  in Theorem \ref{lipaul2}.
These representations exhaust all discrete series representations of $\tilde U(r,s)$ occurring in the discrete spectrum. 
\end{cor}
Generally speaking, in the case $p+q+1=r+s$, there are non-discrete series representations in the discrete spectrum of $\omega|_{\tilde U(r,s)}$.
The complete discrete spectrum can be written down based on \cite{p2}.  
We observe that  $\omega|_{\tilde U(r,s)}$ is no longer $L^2$, but only almost $L^2$. Theorem \ref{l2} does not apply here. Hence there may be limits of the discrete series representations occurring in $ \omega|_{\tilde U(r,s)}^{dis}$.
But we will only deal with the discrete series representations in $\omega|_{\tilde{U}(r,s)}^{dis}$.
\subsection{ $A_{\f q}(\lambda)$ and $\theta(\chi)$}
$A_{\f q}(\lambda)$ is a very important class of representations in the theory of
harmonic analysis and automorphic forms. For example, $A_{\f q}(\lambda)$  contributes to the cohomology of Shimura varieties.
The discrete series representations are all $A_{\f q}(\lambda)$'s with $\f q$ the Borel subalgebra of $\f g_{\mathbb C}$.  We shall refer the readers to the paper by Vogan and Zuckermann \cite{vz} and Vogan \cite{v} for the properties of $A_{\f q}(\lambda)$, and the Knapp-Vogan book \cite{kv} for a systematic construction of these representations. For the application to automorphic forms, we shall refer the readers to Borel-Wallach's book \cite{bw}. In this section, we shall focus on a  class of small $A_{\f q}(\lambda)$ when $\f q$ is maximal parabolic. The branching laws of these $A_{\f q}(\lambda)$ will lead us to a proof of GGP conjecture for $U(p,q)$. \\
\\
We shall follow the notation from \cite{pt}. Let $G=U(r+s, r+s)^{\sqrt{{\det}^{r-s}}}$. Let
$L$ be the corresponding covering of the subgroup $U(r,s) \times U(s,r)$ diagonally embedded in $U(r+s, r+s)$. Fix a standard compact Cartan $\f t$ in $\f u(r,s) \oplus \f u(s,r)$. Fix  a standard positive root system in $\Delta(\f g_{\mathbb C}, i \f t)$, i.e., the positive roots are $\{ e_i- e_j \mid i < j \}$.
Let $\f q$ be the parabolic subalgebra spanned by root vectors from $\f u(r,s)_{\mathbb C}$, $\f u(s,r)_{\mathbb C}$ and positive root vectors. Let $\lambda_{k_1, k_2}$ be the ${\det_{\tilde U(r,s)}}^{k_1} {\det_{\tilde U(s,r)}}^{k_2}$ character of $L$. \underline{Here $k_1$, $k_2$ are integers or half integers and $k_1 \geq k_2$}. By the process known as cohomological induction or Zuckerman's derived functor construction, we obtain $A_{\f q}(\lambda)$ ( see for example \cite{kv}). We denote
this $A_{\f q}(\lambda_{k_1, k_2})$ by $A_{r,s;s,r}(k_1, k_2)$.  \\
\\
Generally speaking, the branching law of $A_{\f q}(\lambda)$ is difficult, due to the lack of a geometric realization.  Even in the discrete series case where $A_{\f q}(\lambda)$ have various geometric realizations (see for example \cite{sch}) and the reference therein), the branching law remains a difficult problem. One important case is the GGP conjecture treated in this paper. Another important case is the branching law of the discrete series restricted to a maximal compact subgroup $K$. The Blattner's formula, proved by Heckt and Schmid, computes the multiplicity of $K$-types in a discrete series. Yet due to the complexity of the Blattner's formula, it is not completely known when the multiplicity is nonzero. Hence even $D(\lambda)|_K$ is not completely known.
One important result about the branching law of $A_{\f q}(\lambda)$ is  due to T. Kobayashi. 
Roughly, Kobayashi determined which of $A_{\f q}(\lambda)$ decomposes discretely when restricted to a  symmetric subgroup (\cite{ko}). In this situation, $A_{\f q}(\lambda)$ decomposes on the Harish-Chandra module level. \\
\\
In general,  $A(r,s;s,r)(k_1, k_2)|_{\tilde U(p,q)}$ will have continuous spectrum. Hence Kobayashi's result does not apply here. Despite of all these difficulties, $A(r,s;s,r)(k_1, k_2)|_{\tilde U(p,q)}$ can be computed by relating it to Howe's correspondence. 
We recall the main result of Paul and Trapa in \cite{pt} (Page 144, Theorem 4.1).
\begin{thm}\label{paultrapa} Let $k$ be an integer.
If $k \geq 0$ then $\theta(r,s; r+s, r+s)({\det}^k) \cong A_{r,s;s,r}(k, 0)$; if $k \leq 0$ then $\theta(r,s; r+s, r+s)({\det}^k) \cong A_{r,s;s,r}(0, k)$.
\end{thm}
We shall remind the reader that $A_{r,s;s,r}(k+n, n)$ differs from $A_{r,s;s,r}(k, 0)$ by a twist of the character ${\det}^n$. Furthermore, by a result of Li (\cite{li}), 
$$\theta(r,s; r+s, r+s)({\det}^k) \cong \omega(r+s, r+s; r,s)^{\infty} \otimes_{U(r,s)^e} {\det}^{-k}.$$
So we obtain
$$A_{r,s;s,r}(k, 0) \cong \omega(r+s, r+s; r,s)^{\infty} \otimes_{U(r,s)^e} {\det}^{-k} \qquad (k \geq 0), $$
$$A_{r,s;s,r}(0, -k) \cong \omega(r+s, r+s; r,s)^{\infty} \otimes_{U(r,s)^e} {\det}^{k} \qquad (k \geq 0).$$
Indeed, the right hand sides are exactly the spaces of smooth vectors.
In the next two sections, we shall derive the branching laws for $A_{r,s;s,r}(k_1, k_2)$ based on these invariant tensor products.

\section{The First Branching Law of $A_{r,s;s,r}(k_1, k_2)$ }
\underline{ We fix $p,q,r,s$}. Throughout this section assume $p+q=r+s$. Let $U(p,q) \times U(q,p)$ be diagonally embedded in $U(p+q, p+q)$. Let $\tilde U(p+q, p+q)$ be $U(p+q,p+q)^{\sqrt{{\det}^{r-s}}}$. So $\tilde{U}(p,q)$ will be just $U(p,q)^{\sqrt{{\det}^{r-s}}}$.  We consider $A_{r,s;s,r}(k_1, k_2)|_{ \tilde U(p,q)}$. In this section, we may simply write $A(k_1, k_2)$ for $A_{r,s;s,r}(k_1, k_2)$.  By the remark following Theorem \ref{paultrapa}, it suffices to discuss the branching law for $A(0, -k) (k \geq 0)$. We know that $A(k_1, k_2)|_{\tilde U(p,q)}$ are square integrable, by employing the matrix coefficients estimate based on the Langland parameter. By Theorem \ref{l2}, only discrete series can appear in the discrete spectrum of $A(k_1, k_2)|_{\tilde U(p,q)}$.  Therefore discrete spectrum of $A(k_1, k_2)|_{ \tilde U(p,q)}$ can be obtained by studying $\mc H_{\sigma^c}^{\infty} \otimes_{\tilde U(p,q)} A(k_1, k_2)^{\infty}$ for each discrete series representation $\sigma$.\\
\\
Suppose $k \geq 0$. Notice that $ \omega(r+s, r+s; r,s)^{\infty} \otimes_{U(r,s)^e} {\det}^{k}$ is the space of smooth vectors in $ A(0, -k)$. It suffices to compute 
$$\mc H_{\sigma^c}^{\infty} \otimes_{\tilde U(p,q)} [\omega(r+s, r+s; r,s)^{\infty} \otimes_{U(r,s)^e} {\det}^{k}].$$
In order to apply the associativity of invariant tensor products (Theorem \ref{asso}),  we must give an estimate about the decaying of matrix coefficients of $\omega(r+s, r+s; r,s)^{\infty}|_{\tilde U(p,q) \times \tilde U(r,s)}$. We can then obtain a direct estimate on certain smooth matrix coefficients of $A(0,k)$ which will be weaker than the estimate from the Langlands parameters, but sufficient for the purpose of this paper. Consequently, $\mc H_{\sigma^c}^{\infty} \otimes_{\tilde U(p,q)} A(k_1, k_2)^{\infty}$ can be expressed as a composition of Howe's correspondences, which can then be computed in terms of the Harish-Chandra parameter. The ideas here are similar to \cite{basic}.
\subsection{Matrix coefficient Estimates: General Case}
Let $G$ be a  real reductive group with compact center. Fix a Cartan involution $\Theta$ on $G$. Let $K$ be the maximal compact subgroup fixed by $\Theta$. Denote the Cartan involution on $\f g$ also by $\Theta$.  Let $\f p$ be the $-1$ eigenspace of $\Theta$. Then
$\f g=\f k \oplus \f p$. Fix a maximal commutative subalgebra $\f a$ of $\f p$. Let $A$ be the the analytic subgroup generated by $\f a$. Fix a Weyl chamber $A^+$ in $A$. 
Then $G$ has a Cartan decomposition of the form $KA^+K$. For the noncompact unitary groups, $A^+$ component of $g$ will be unique. \underline{Write the $A^+$-component of $g \in G$ as $\exp H(g)$ with $H(g) \in \f a^+$}. The invariant measure of $G$ can be expressed as $\Delta(H) dK d H d K$. The function $\Delta(H)$ has the exponential growth with exponent
$2 \rho(G)$. Here $2 \rho(G)$ is the sum of all positive restricted roots corresponding to the Weyl chamber $\f a^+$. For a function $f$ to be integrable, it suffices that $f(g)$  \lq\lq decays" faster than $\exp - 2\rho(G)(H(g))$. For a  a function $f$ to be square integrable, it suffices that $f(g)$  \lq\lq decays" faster than $\exp - \rho(G)(H(g))$.\\
\\
Consider now the oscillator representation $\omega$ of $\widetilde{Sp}_{2N}(\mb R)$. Choose 
$$\f a= \{ H=\diag(H_1, H_2, \ldots H_N, -H_1, -H_2, \ldots, -H_N) \mid H_i \in \mathbb R \}$$ and $\f a^+=\{ H \mid H_1 \geq H_2 \geq H_3 \ldots \geq H_N \geq 0 \}$. Then
$A^+$ consists of 
$$\diag(a_1, a_2, \ldots a_N, a_1^{-1}, a_2^{-1}, \ldots a_N^{-1}) \qquad a_1 \geq a_2 \geq \ldots a_N \geq 1.$$
\uline{For simplicity, we write $g \cong \diag(a_1^{\pm 1}, a_2^{\pm 1}, \dots a_N^{\pm 1})$ if $g=
k \diag(a_1^{\pm 1}, a_2^{\pm 1}, \dots a_N^{\pm 1}) k^{\prime}$ with $k, k^{\prime} \in K$}. \\
\\
It is well-known that the smooth matrix coefficients of $\omega$ satisfy
$$f(g) \leq C \prod_{i=1}^N (a_i+a_i^{-1})^{-\frac{1}{2}}.$$
See for example the proof of Theorem 3.3 \cite{theta} Pg. 264. Here $C$ is a constant that depends on $f$.
Now the dual reductive pair $(U(m,n), U(r,s))$ is embedded in $Sp_{2(m+n)(r+s)}(\mb R)$. Suppose that $m \geq n$ and $ r \geq s$.  Let $g_1 \in \tilde U(m,n)$ and $g_2 \in \tilde U(r, s)$. 
Let 
$$g_1 \cong \diag(a_1^{\pm 1}, a_2^{\pm 1}, \ldots a_n^{\pm 1}, \overbrace{1, \ldots, 1}^{m-n}), \qquad
g_2 \cong \diag(b_1^{\pm 1}, b_2^{\pm 1}, \ldots b_{s}^{\pm 1}, \overbrace{1, \ldots, 1}^{r-s}).$$
In $\widetilde{Sp}_{2(m+n)(r+s)}(\mb R)$,
$$g_1 g_2 \cong \diag(\{ a_i^{\pm 1} b_j^{\pm 1} \}_{i \in [1,n], j \in [1,s]}, \{ a_i^{\pm 1} b_j^{\pm 1} \}_{i \in [1,n], j \in [1,s]}, \overbrace{ \{ a_i^{\pm 1} \}_{i \in [1,n]}, \ldots, \{ a_i^{\pm 1} \}_{i \in [1,n]} }^{2r-2s}, $$
$$ \overbrace{ \{ b_j^{\pm 1} \}_{j \in [1,s]}, \ldots, \{ b_j^{\pm 1} \}_{j \in [1,s]}}^{2m-2n}, \overbrace{1, \ldots, 1}^{2(m-n)(r-s)} ).$$
Then smooth matrix coefficients of $\omega(m,n;r,s)$ restricted onto $\tilde U(m,n) \times \tilde U(r, s)$ is bounded by
$$|f(g_1 g_2)| \leq C [\prod_{i=1}^n \prod_{j=1}^s (a_i b_j+a_i^{-1} b_j^{-1})^{-1} (a_i b_j^{-1}+a_i^{-1} b_j)^{-1}] [\prod_{i=1}^n (a_i+ a_i^{-1})^{-1} ]^{r-s} [\prod_{j=1}^s (b_j+b_j^{-1})^{-1}]^{m-n}.$$
We summarize our result in the following lemma.
\begin{lem}
Suppose $m \geq n$ and $r \geq s$. Let $f$ be a smooth matrix coefficient of $\omega(m,n;r,s)$. Let $g_1 \in \tilde U(m,n)$ and $g_2 \in \tilde U(r,s)$. Then
$$ f(g_1 g_2) \leq C [\prod_{i=1}^n \prod_{j=1}^s (a_i^2+ b_j^2+a_i^{-2}+ b_j^{-2})^{-1} ] [\prod_{i=1}^n (a_i+ a_i^{-1}) ]^{-r+s} [\prod_{j=1}^s (b_j+b_j^{-1})]^{-m+n}.$$
\end{lem}
\subsection{Matrix Coefficients of $\omega(n, n; r,s)$}
Let $m=n=r+s$. Let $p+q=n$ and $U(p,q) \times U(q,p)$ be diagonally embedded in $U(n,n)$. Let $g_1 \in \tilde U(p,q)$. For the purpose of matrix coefficients estimation, we assume $p \geq q$ and $r \geq s$. This will cause no loss of generality. Then as an element in $U(n,n)$,
$$g_1 \cong \diag(a_1^{\pm 1}, a_2^{\pm 1}, \ldots a_q^{\pm 1}, \overbrace{1, 1, \ldots 1}^{2p}).$$
We have for any $g_2 \in \tilde U(r,s)$,
\begin{equation}\label{f(g_1g_2)}
f(g_1 g_2) \leq C \left[ \prod_{i=1}^q \prod_{j=1}^s (a_i^2+ b_j^2+a_i^{-2}+ b_j^{-2})^{-1} \right]
\left[\prod_{j=1}^s (b_j^2+b_j^{-2}+2) \right]^{-p} \left[\prod_{i=1}^q (a_i+ a_i^{-1}) \right]^{-r+s} .
\end{equation}
When $q=0$, $g_1$ is in a compact group. We have
$$f(g_1 g_2) \leq C \left[\prod_{j=1}^s (b_j^2+b_j^{-2}+2) \right]^{-n} .$$

\begin{thm}\label{est} Suppose $p+q=r+s=n$. Let $f$ be a matrix coefficient of smooth vectors for $\omega(n,n; r,s)$. Let $g_1 \in \tilde U(p,q) \subseteq \tilde U(n,n)$ and $g_2 \in \tilde U(r,s)$. Let $v(g_1)$ be a $L^2$-function on $\tilde U(p,q)$. Then $ |v(g_1) f(g_1 g_2)| \in L^1(\tilde U(p,q) \times \tilde U(r,s))$.
\end{thm}
Proof: Without loss of generality, assume $p \geq q$ and $r \geq s$. This will cause no problem with the estimates of matrix coefficients. Let 
$g_1 \cong \diag(a_1^{\pm 1}, a_2^{\pm 1}, \ldots, a_q^{\pm 1}, \overbrace{1, \ldots 1}^{p-q})$ with $a_1 \geq a_2 \geq \ldots \geq a_q \geq 1$. Notice 
$\Delta(g_1) \leq C \Pi_{i=1}^q a_i^{2p+2q-4 i+2}.$ The exponents here come from
$$ 2 \rho(\tilde U(p,q))=(2p+2q-2, 2p+2q-6, \ldots, 2p-2q+2).$$
Let $g_2 \cong \diag(b_1^{\pm 1}, b_2^{\pm 1}, \ldots, b_{s}^{\pm 1}, \overbrace{1, \ldots, 1}^{r-s})$
with $b_1 \geq b_2 \geq \ldots \geq b_s \geq 1$. Notice that
$\Delta(g_2) \leq C \Pi_{j=1}^s b_j^{2r+2s-4 j+2}.$ The exponents here come from
$$ 2 \rho(\tilde U(r,s))=(2r+2s-2, 2r+2s-6, \ldots, 2r-2s+2).$$
Observe that for any $\alpha \in [0,1]$,
$$(a_i^2+a_i^{-2}+b_j^2+b_j^{-2})^{-1}  \leq (a_i+ a_i^{-1})^{- 2 \alpha} (b_j+ b_j^{-1})^{-2+ 2 \alpha}.$$
We have three cases.
\begin{enumerate}

\item Suppose that $q \geq s \geq 1$. Then we have
$$ \prod_{i=1}^q (a_i^2+b_1^2+a_i^{-2}+ b_1^{-2})^{-1} \leq (a_1+a_1^{-1})^{-2+\frac{1}{2}}
(b_1+b_1^{-1})^{-2q+2 -\frac{1}{2}}, $$
$$ \prod_{i=1}^q (a_i^2+b_2^2+a_i^{-2}+ b_2^{-2})^{-1} \leq (a_2+a_2^{-1})^{-2+\frac{1}{2}}(a_1+a_1^{-1})^{-2}
(b_2+b_2^{-1})^{-2q+4 -\frac{1}{2}} ,$$
$$\ldots \ldots \qquad \ldots \ldots \qquad \ldots \ldots $$
$$\prod_{i=1}^q (a_i^2+b_s^2+a_i^{-2}+ b_s^{-2})^{-1} \leq (a_s+a_s^{-1})^{-2+\frac{1}{2}}(a_{s-1}+a_{s-1}^{-1})^{-2} \ldots (a_1+a_1^{-1})^{-2}
(b_s+b_s^{-1})^{-2q+2s -\frac{1}{2}} .$$
Multiplying them, we obtain
$$ \prod_{i=1}^q \prod_{j=1}^s (a_i^2+ b_j^2+a_i^{-2}+ b_j^{-2})^{-1}  \leq \prod_{i=1}^s (a_i+a_i^{-1})^{-2s+2i-2+\frac{1}{2}} \prod_{j=1}^s (b_j+b_j^{-1})^{-2q+2j -\frac{1}{2}}.$$
Combined with Equation \ref{f(g_1g_2)}, we have
$$|f(g_1 g_2)| \leq C \prod_{i=1}^s (a_i+a_i^{-1})^{2i-2+\frac{1}{2}-r-s} \prod_{j=1}^s (b_j+b_j^{-1})^{-2p-2q+2j -\frac{1}{2}} \prod_{i=s+1}^q (a_i+a_i^{-1})^{-r+s}.$$
Since $p+q=r+s$ and
$$\prod_{j=1}^s (b_j+b_j^{-1})^{-2r-2s+2j -\frac{1}{2}} \in L^1(\tilde U(r,s), \Delta(b) d b \, d k_1 \, d k_1^{\prime}), $$
$$ 
\prod_{i=1}^s (a_i+a_i^{-1})^{2i-2+\frac{1}{2}-p-q} \prod_{i=s+1}^q (a_i+a_i^{-1})^{-r+s} \in L^2(\tilde U(p,q), \Delta(a) d a \, d k_2 \, d k^{\prime}_2),$$
it follows that $ |v(g_1) f(g_1 g_2)| \in L^1(\tilde U(p,q) \times \tilde U(r,s))$.
\item Suppose that $s \geq q \geq 1$. Similarly,  we  have
$$\prod_{j=1}^s (a_1^2+b_j^2+a_1^{-2} + b_j^{-2})^{-1} \leq (a_1+a_1^{-1})^{-2s+\frac{1}{2}} (b_1+ b_1^{-1})^{-\frac{1}{2}} , $$
$$\prod_{j=1}^s (a_2^2+b_j^2+a_2^{-2} + b_j^{-2})^{-1} \leq (a_2+a_2^{-1})^{-2s+2+\frac{1}{2}} (b_{2}+ b_{2}^{-1})^{-\frac{1}{2}} (b_1+ b_1^{-1})^{-2} ,$$
$$\ldots \ldots \qquad \ldots \ldots \qquad \ldots \ldots $$
$$\prod_{j=1}^s (a_q^2+b_j^2+a_q^{-2} + b_j^{-2})^{-1} \leq (a_q+a_q^{-1})^{-2s+2q-2+\frac{1}{2}} (b_q+ b_q^{-1})^{-\frac{1}{2}}(b_{q-1}+b_{q-1}^{-1})^{-2} \ldots (b_1+b_1^{-1})^{-2}.  $$
Multiplying them, we obtain
$$\prod_{i=1}^q \prod_{j=1}^s (a_i^2+b_j^2+a_i^{-2} + b_j^{-2})^{-1}  \leq \prod_{i=1}^q
 (a_i+a_i^{-1})^{-2s+2i-2+\frac{1}{2}} \prod_{j=1}^q (b_j+ b_j^{-1})^{-2q+2j-\frac{1}{2}}. $$
Combined with Equation \ref{f(g_1g_2)}, we have
$$|f(g_1 g_2)| \leq C \prod_{i=1}^q (a_i+a_i^{-1})^{2i-2+\frac{1}{2}-r-s} \prod_{j=1}^q (b_j+b_j^{-1})^{-2p-2q+2j -\frac{1}{2}} \prod_{j=q+1}^s (b_j+b_j^{-1})^{-2p}.$$
Since $p+q=r+s$ and
$$\prod_{i=1}^q (a_i+a_i^{-1})^{2i-2+\frac{1}{2}-p-q} \in L^2(\tilde U(p,q), \Delta(a) d a \,  d k_2 \, d k^{\prime}_2),$$
$$\prod_{j=1}^q (b_j+b_j^{-1})^{-2r-2s+2j -\frac{1}{2}} \prod_{j=q+1}^s (b_j+b_j^{-1})^{-2p} \in  L^1(\tilde U(r,s), \Delta(b) d b \, d k_1 \, d k_1^{\prime}), $$
we see that $ |v(g_1) f(g_1 g_2)| \in L^1(\tilde U(p,q) \times \tilde U(r,s))$.
\item If $s=0$ or $q=0$, our assertion becomes obvious by Equation \ref{f(g_1g_2)}.
\end{enumerate}
$\Box$\\
\\
Our proof shows for vectors in the smooth representation $\omega^{\infty}(n,n;r,s) \otimes_{\tilde U(r,s)} \det^{k}$ of $\tilde U(n,n)$, the matrix coefficients restricted to $\tilde U(p,q)$ will be bounded by
$$C \prod_{i=1}^q (a_i+a_i^{-1})^{2i-2+\frac{1}{2}-p-q} \in L^2(\tilde U(p,q), \Delta(a) d a \,  d k_2 \, d k^{\prime}_2).$$
\begin{cor}\label{estimate} If $p+q=r+s$, both $A_{r,s;s,r}(k_1, k_2)|_{\tilde U(p,q)}$ and  $A_{r,s;s,r}(k_1, k_2)|_{\tilde U(p-1,q)} (p \geq 1)$ are square integrable representations.
\end{cor}

\subsection{Associativity and discrete spectrum of $A(0, -k)$}
Let $\sigma$ be a discrete series representation of $\tilde U(p,q)$.  By Lemma \ref{basicproperty}, Theorem \ref{asso} and Theorem \ref{multi}, we compute the multiplicity space 
\begin{equation}
\begin{split}
 & M_{A(0, -k)}(\sigma) \\
= & \mc H_{\sigma^c}^{\infty} \otimes_{\tilde U(p,q)} [ \omega(p+q, p+q; r,s)^{\infty} \otimes_{\tilde U(r,s)} {\det}^{k} ] \\
\cong &  \mc H_{\sigma^c}^{\infty} \otimes_{\tilde U(p,q)} [ (\omega(p,q;r,s)^{\infty} \otimes \omega(q,p;r,s)^{\infty}) \otimes_{\tilde U(r,s)} {\det}^{k}] \\
 \cong & [ \mc H_{\sigma^c}^{\infty} \otimes_{\tilde U(p,q)} ( \omega(p,q;r,s)^{\infty}  \otimes \omega(q,p;r,s)^{\infty}) ] \otimes_{\tilde U(r,s)} {\det}^k  \\
 \cong & [ (\mc H_{\sigma^c}^{\infty} \otimes_{\tilde U(p,q)}  \omega(p,q;r,s)^{\infty} ) \otimes \omega(q,p;r,s)^{\infty} ] \otimes_{\tilde U(r,s)} {\det}^k  \\
 \cong  & [ (\mc H_{\sigma^c}^{\infty} \otimes_{\tilde U(p,q)}  \omega(p,q;r,s)^{\infty} ) \otimes {\det}^k ] \otimes_{\tilde U(r,s)} \omega(q,p;r,s)^{\infty} \\
 \cong & [\theta(p,q;r,s)(\sigma)^{\infty} \otimes {\det}^k ] \otimes_{\tilde U(r,s)} \omega(q,p;r,s)^{\infty} \\
 \cong & \theta(r,s;q,p)([\theta(p,q;r,s)(\sigma) \otimes {\det}^k]^c).
\end{split}
\end{equation}
The last two equivalences follow from  Theorem \ref{lipaul}. All the equivalences shall be interpreted as equivalence of the completions as irreducible unitary representations, following the convention set in 1.1. By Theorem \ref{multi}, we have
\begin{lem}\label{atheta} $\sigma \in A_{r,s;s,r}(0,-k)|_{\tilde U(p,q)}$ if and only if
$\theta(r,s;q,p)([\theta(p,q;r,s)(\sigma) \otimes {\det}^k]^*) \neq 0$.
\end{lem}
Recall that the Harish-Chandra parameters for $U(p,q)^e$ ($p+q$ even) and $U(p,q)^{o}$ ($p+q$ odd) are half integers.  Assume that $\theta(p,q;r,s)(\sigma)=D(\lambda^+, \lambda^-)$. Then $\sigma=D((\lambda_{+}^{+}, \lambda_{-}^{-}), (\lambda^{-}_{+}, \lambda^{+}_{-}))$ and we must have 
$$\ca(\lambda^+)=r, \qquad \ca(\lambda^-)=s, \qquad \ca(\lambda^{+}_+)+\ca(\lambda_{-}^{-})=p, \qquad \ca(\lambda^{+}_-)+\ca(\lambda_{+}^{-})=q.$$
Clearly, 
$\theta(p,q;r,s)(\sigma) \otimes {\det}^k =D(\lambda^+ + \mathbf{k}, \lambda^- +\mathbf{k})$. Hence
$$ [\theta(p,q;r,s)(\sigma) \otimes {\det}^k]^c= D(-(\lambda^+ + \mathbf{k}), -(\lambda^- +\mathbf{k})).$$
By Theorem \ref{lipaul},
$\theta(r,s;q,p)([\theta(p,q;r,s)(\sigma) \otimes {\det}^k]^c) \neq 0$ if and only if 
$$\ca((\lambda^+ +\mathbf{k})_{-})+ \ca((\lambda^{-}+ \mathbf{k})_{+})=q, \qquad
\ca((\lambda^+ +\mathbf{k})_{+})+ \ca((\lambda^{-}+ \mathbf{k})_{-})=p.$$
Comparing with $ \ca(\lambda^{+}_-)+\ca(\lambda_{+}^{-})=q$ and $\ca(\lambda^{+}_+)+\ca(\lambda_{-}^{-})=p$, it is sufficient and necessary that
$$\ca(\lambda^{+} \cap [-k, 0])=\ca(\lambda^{-} \cap [-k, 0]).$$
Then we have
$$\theta(r,s;q,p)([\theta(p,q;r,s)(\sigma) \otimes {\det}^k]^c)= D((-([\lambda^++ \mathbf{k}]_{-}), -([\lambda^{-}+ \mathbf{k}]_{+})), (-([\lambda^{-}+\mathbf{k}]_{-}), -([\lambda^++ \mathbf{k}]_{+}))).$$
We have proved
\begin{thm}\label{discrete:equalcase} Suppose that $r+s=p+q$ and $k \geq 0$. The discrete spectrum
$$A_{r,s;s,r}(0, -k)|_{\tilde U(p,q)}^{dis} = A_{r,s;s,r}(0, -k)|_{\tilde U(p, q) \times \tilde U(q,p)}^{dis} = A_{r,s;s,r}(0, -k)|_{\tilde U(q,p)}^{dis} $$
is the direct sum of
$$D((\lambda_{+}^{+}, \lambda_{-}^{-}), (\lambda^{-}_{+}, \lambda^{+}_{-})) \boxtimes D(([\lambda^++ \mathbf{k}]_{-}, [\lambda^{-}+ \mathbf{k}]_{+}), ([\lambda^{-}+k]_{-}, [\lambda^++ \mathbf{k}]_{+}))^*$$
where the entries of $(\lambda^+, \lambda^-)$ are all half integers and satisfy the following equations
$$\ca(\lambda^+)=r, \qquad \ca(\lambda^-)=s, \qquad \ca(\lambda^{+}_+)+\ca(\lambda_{-}^{-})=p, \qquad \ca(\lambda^{+}_-)+\ca(\lambda_{+}^{-})=q,$$
$$\ca(\lambda^{+} \cap [-k, 0])=\ca(\lambda^{-} \cap [-k, 0]).$$
\end{thm}
As to the continuous spectrum, the situation is similar to the symmetric space case. That is, the continuous spectrum will comes from induced representations of the discrete series that appear in the discrete spectrum of a smaller $A(0, -k)$.  A somewhat easier approach is to prove Lemma \ref{atheta} for $\sigma \in_{wk} A_{r,s;s,r}(0,-k)|_{\tilde U(p,q)}$. Recall that the continuous spectrum is only defined almost everywhere. Therefore, up to a measure zero set, one can obtain the complete spectrum of $A_{r,s;s,r}(0,-k)|_{\tilde U(p,q) \times \tilde U(q,p)}$ this way. However, the whole support of $A_{r,s;s,r}(0,-k)|_{\tilde U(p,q) \times \tilde U(q,p)}$ is more delicate and requires more caution. We will not discuss it in this paper.

\section{The Second Branching Law }
\underline{We still fix $p+q=r+s$.  Assume $p \geq 1$ and $k \geq 1$}. Let $U(p-1,q) \times  U(q+1, p)$ be diagonally embedded in $U(p+q, p+q)$. We consider $A(0, -k)|_{\tilde U(p-1,q) \times \tilde U(q+1, p)}$. By Cor \ref{estimate}, $A(0,-k)|_{\tilde U(p-1,q)}$ is square integrable. Thus only discrete series can occur in $A(0,-k)|_{\tilde U(p-1,q)}^{dis}$. Let $\sigma$ be a discrete series representation of $\tilde U(p-1,q)$. The entries of the Harish-Chandra parameters involved here will all be integers. Recall from Theorem \ref{multi}, the multiplicity space
$$M_{A(0, -k)}(\sigma) \cong \mc H_{\sigma^c}^{\infty} \otimes_{\tilde U(p-1,q)}  A(0,-k)^{\infty}.$$
Clearly
$$\mc H_{\sigma^c}^{\infty} \otimes_{\tilde U(p-1,q)}  A(0,-k)^{\infty} \cong \mc H_{\sigma^c}^{\infty} \otimes_{\tilde U(p-1,q)} [(\omega(p-1,q; r,s)^{\infty} \otimes \omega(q+1, p; r,s)^{\infty}) \otimes_{\tilde U(r,s)} {\det}^k].$$
By essentially the same estimate in Theorem \ref{est}, we can apply the law of associativity:
$$\mc H_{\sigma^c}^{\infty} \otimes_{\tilde U(p-1,q)}  A(0,-k)^{\infty} \cong [(\mc H_{\sigma^c}^{\infty} \otimes_{\tilde U(p-1,q)} \omega(p-1,q; r,s)^{\infty}) \otimes \omega(q+1, p; r,s)^{\infty}] \otimes_{\tilde U(r,s)} {\det}^k.$$
By associativity and commutativity, we have
\begin{equation}
\begin{split}
 \mc H_{\sigma^c}^{\infty} \otimes_{\tilde U(p-1,q)}  A(0,-k)^{\infty} \cong & [(\mc H_{\sigma^c}^{\infty} \otimes_{\tilde U(p-1,q)} \omega(p-1,q; r,s)^{\infty} ) \otimes {\det}^k ] \otimes_{\tilde U(r,s)} \omega(q+1, p; r,s)^{\infty}.\\
\cong & [\theta(p-1,q; r,s)(\sigma)^{\infty} \otimes {\det}^k] \otimes_{\tilde U(r,s)} \omega(q+1, p; r,s)^{\infty} \\
\cong & \theta(r,s; q+1,p)( [\theta(p-1,q; r,s)(\sigma) \otimes {\det}^k]^c).
\end{split}
\end{equation}
For technical reasons, we need to assume that the Harish-Chandra parameter of $\sigma$ does not have any zeros. Then $\mc H_{\sigma^c}^{\infty} \otimes_{\tilde U(p-1,q)} A(k_1, k_2)^{\infty}$ is a discrete series of $\tilde{U}(q+1,p)$.  Our equations are valid due to Theorem \ref{lipaul2}.

\begin{prop}\label{5.1}
Suppose that $\sigma$ is a discrete series representation of $\tilde U(p-1, q)$ and its Harish-Chandra parameter does not contain any zeros.
 $\sigma \in A(0, -k)|_{\tilde U(p-1,q)}$ if and only if   the multiplicity space
$$M_{A(0,-k)}(\sigma) \cong \theta(r,s; q+1,p)( [\theta(p-1,q; r,s)(\sigma) \otimes {\det}^k]^c) \neq 0.$$
In this situation, this multiplicity space is an irreducible tempered representation of $\tilde U(q+1, p)$. 
\end{prop}
Proof: The irreducibility of $M_{A(0,-k)}(\sigma)$ follows from Howe's theorem. The temperedness of $M_{A(0,-k)}(\sigma)$ follows from a theorem of Paul (Prop. 1.4 \cite{p2}). This also proves that $A(0,-k)|^{dis}_{\tilde U(p-1,q) \times \tilde U(q+1, p)}$  has multiplicity one. $\Box$ \\
\\
This proposition gives us the complete description of the discrete spectrum of $ A(0, -k)|_{\tilde U(p-1,q) \times \tilde U(q+1, p)}$ in terms of Howe's correspondence. Similar statements are true for all other subgroups that are diagonally embedded in $U(n,n)$. Since the main focus here is GGP conjecture and the discrete series, we assume: \\
\\
\uline{the multiplicity space $M_{A(0, -k)}(\sigma)$ is equivalent to a discrete series representation of $\tilde U(q+1,p)$}. \\
\\
This assumption excludes those $\sigma$ in the discrete series of $\tilde U(p-1,q)$ with $M_{A(0,-k)}(\sigma)$ a limit discrete series of $\tilde U(q+1,p)$. \\
\\
It follows form Theorem \ref{lipaul2} that $\theta(p-1,q; r,s)(\sigma)$ is a discrete series and 
\begin{enumerate}
\item either
 $\theta(p-1,q; r,s)(\sigma)= D((\mu^+, 0), \mu^{-}) $
 \item  or $ \theta(p-1,q; r,s)(\sigma)= D(\mu^+, (\mu^{-},0)) .$
 \end{enumerate}
 In both cases $\sigma= D((\mu_{+}^{+}, \mu_{-}^{-}), (\mu^{-}_{+}, \mu_{-}^{+}))$ and all entries of $\mu_{\pm}^{\pm}$ are nonzero. But the cardinalities of $\mu^{\pm}_{\pm}$ will be different. These two cases are mutually exclusive.

 \subsection{Computation on Harish-Chandra Parameters: Case I}
 Suppose that $\theta(p-1,q; r,s)(\sigma)= D((\mu^+, 0), \mu^{-}) $. Then
 $$[\theta(p-1,q; r,s)(\sigma) \otimes {\det}^k]^c = D((-\mu^+ - {\bf{k}}, - k), (-\mu^{-}-{\bf{k}})).$$
 
 \begin{lem}\label{branching21} Let $\sigma= D((\mu_{+}^{+}, \mu_{-}^{-}), (\mu^{-}_{+}, \mu_{-}^{+}))$ with $\mu^{\pm}_{\pm}$ all integers. Assume that 
 $$\ca(\mu_{+}^{+})+ \ca(\mu_{-}^+)+1=r, \qquad \ca(\mu_{+}^{-})+\ca(\mu_{-}^{-})=s, \qquad \ca(\mu_{+}^{+})+ \ca(\mu_{-}^-)=p-1, \qquad \ca(\mu_{+}^{-})+\ca(\mu_{-}^{+})=q$$
 and none of the entries of $\mu$ equals  $0$ or $-k$.
Then $M_{A(0,-k)}(\sigma)$ equals
 \begin{enumerate} 
 \item 
 $$D((-[(\mu^+ + {\bf{k}},  k)_{-}],  0, -[(\mu^{-}+{\bf{k}})_{+}]), (-[(\mu^{-}+{\bf{k}})_{-}], -[(\mu^+ + {\bf{k}},     k)_{+}]))$$
 when $\ca(\mu^{+} \cap (-k, 0)) = \ca(\mu^{-} \cap (-k, 0))$;
 \item 
 $$D((-[(\mu^+ + {\bf{k}},  k)_{-}], -[(\mu^{-}+{\bf{k}})_{+}]), (-[(\mu^{-}+{\bf{k}})_{-}],0, -[(\mu^+ + {\bf{k}},     k)_{+}]))$$
 when $1+ \ca(\mu^{+} \cap (-k, 0)) = \ca(\mu^{-} \cap (-k, 0))$; 
 \item zero otherwise.
 \end{enumerate}
 \end{lem}
 Proof: By Theorem \ref{lipaul2} and Prop. \ref{5.1}, (1) holds iff 
 $$\ca((\mu^+ + {\bf{k}},  k)_{-})+\ca((\mu^{-}+{\bf{k}})_{+}) =q, \qquad
 \ca((\mu^{-}+{\bf{k}})_{-})+\ca((\mu^+ + {\bf{k}},     k)_{+}) = p. $$
 By comparing with 
 $\ca(\mu_{+}^{+})+ \ca(\mu_{-}^-)=p-1$ and $\ca(\mu_{+}^{-})+\ca(\mu_{-}^{+})=q,$
 (1) holds iff $\ca(\mu^{+} \cap (-k, 0)) = \ca(\mu^{-} \cap (-k, 0))$. \\
 \\
 Similarly, (2) holds iff 
 $$\ca((\mu^+ + {\bf{k}},  k)_{-})+\ca((\mu^{-}+{\bf{k}})_{+}) =q+1, \qquad
 \ca((\mu^{-}+{\bf{k}})_{-})+\ca((\mu^+ + {\bf{k}},  k)_{+}) = p-1. $$ 
  if and only if $1+\ca(\mu^{+} \cap (-k, 0)) = \ca(\mu^{-} \cap (-k, 0))$. \\
 \\
 By Theorem \ref{lipaul2}, (3) also holds. $\Box$
 
 \subsection{Computation on Harish-Chandra Parameters: Case II}
 Suppose that $\theta(p-1,q; r,s)(\sigma)= D(\mu^+, (\mu^{-}, 0)) $. Then
 $$[\theta(p-1,q; r,s)(\sigma) \otimes {\det}^k]^c = D(-(\mu^+ + {\bf{k}}), -(\mu^{-}+{\bf{k}},k)).$$
 \begin{lem}\label{branching22} Let $\sigma= D((\mu_{+}^{+}, \mu_{-}^{-}), (\mu^{-}_{+}, \mu_{-}^{+}))$ with $\mu_{\pm}^{\pm}$ all integers. Assume that 
 $$\ca(\mu_{+}^{+})+ \ca(\mu_{-}^+)=r, \qquad \ca(\mu_{+}^{-})+\ca(\mu_{-}^{-})+1=s, \qquad \ca(\mu_{+}^{+})+ \ca(\mu_{-}^-)=p-1, \qquad \ca(\mu_{+}^{-})+\ca(\mu_{-}^{+})=q$$
 and none of the entries of $\mu$ equals $0$ or $-k$.
 We have $M_{A(0, -k)}(\sigma)$ equals
 \begin{enumerate} 
 \item 
 $$D((-[(\mu^+ + {\bf{k}})_{-}],  0, -[(\mu^{-}+{\bf{k}},k)_{+}]), (-[(\mu^{-}+{\bf{k}}, k)_{-}], -[(\mu^+ + {\bf{k}})_{+}]))$$
 when $\ca(\mu^{+} \cap (-k, 0)) = 1+\ca(\mu^{-} \cap (-k, 0))$;
 \item 
 $$D((-[(\mu^+ + {\bf{k}})_{-}], -[(\mu^{-}+{\bf{k}},k)_{+}]), (-[(\mu^{-}+{\bf{k}},k)_{-}],0, -[(\mu^+ + {\bf{k}})_{+}]))$$
 when $ \ca(\mu^{+} \cap (-k, 0)) = \ca(\mu^{-} \cap (-k, 0))$; 
 \item zero otherwise.
 \end{enumerate}
 \end{lem}
\subsection{Discrete Spectrum of $A(0, -k)|_{\tilde U(p-1, q)}$}
Let $\sigma$ be a discrete series representation of $\tilde U(p-1, q)$. Recall that $\sigma \boxtimes M_{A(0,-k)}(\sigma)$ will be in the discrete spectrum of $A(0,-k)|_{\tilde U(p-1,q)}$ if and only if $M_{A(0, -k)}(\sigma)$ is nonvanishing. 
\begin{thm}\label{branching2} Let $\sigma \boxtimes M_{A(0,-k)}(\sigma)$ be in the discrete spectrum of $A_{r,s;s,r}(0,-k)|_{\tilde U(p-1,q) \times \tilde U(q+1,p)}$. Suppose that $M_{A(0,-k)}(\sigma)$ is a discrete series representation of $\tilde U(q+1, p)$. Then $\sigma$ and $M_{A(0,-k)}(\sigma)$ must be given by the recipe in Lemma \ref{branching21} or Lemma \ref{branching22}.
\end{thm}
Proof: By Prop. \ref{5.1}, $M_{A(0,-k)}(\sigma)=\theta(r,s;q+1,p)( [\theta(p-1,q; r,s)(\sigma) \otimes {\det}^k]^c)$. By Theorem \ref{lipaul2}, since $\theta(r,s;q+1,p)( [\theta(p-1,q; r,s)(\sigma) \otimes {\det}^k]^c)$ is a discrete series representation, $\theta(p-1,q; r,s)(\sigma)$ must be a discrete series representation and $\sigma$ must also be a discrete series representation. In addition, the Harish-Chandra parameters of $\sigma$ and $\theta(p-1,q; r,s)(\sigma) \otimes {\det}^k$ cannot contain any zero.  Let 
$\sigma=D((\mu_{+}^{+}, \mu_{-}^{-}), (\mu^{-}_{+}, \mu_{-}^{+}))$. Then $\theta(p-1,q; r,s)(\sigma) \otimes {\det}^k= D(\mu^+ +\mathbf{k}, (\mu^{-}
+\mathbf{k},k))$ or $D((\mu^+ +\mathbf{k}, k), \mu^{-} +\mathbf{k})$. Hence none of the entries of $\mu^{\pm}$ can be equal to $-k$.  The equalities in Lemma \ref{branching21} and Lemma \ref{branching22} are guaranteed by $\theta(r,s;q+1,p)( [\theta(p-1,q; r,s)(\sigma) \otimes {\det}^k]^c) \neq 0$. Our assertion then follows. $\Box$ \\
\\
Now we are interested in some special discrete series representation $M_{A(0,-k)}(\sigma)$,
namely those whose Harish-Chandra parameters do not contain any integers in $(-k, 0)$. By Prop. \ref{5.1} and Theorem \ref{lipaul2}, we must have 
$$\ca(\mu^{+} \cap [-k, 0]) = \ca(\mu^{-} \cap [-k, 0])=0.$$
Combining with Lemma \ref{branching21} and Lemma \ref{branching22}, we obtain
\begin{cor}\label{induction1} Suppose that $M_{A(0,-k)}(\sigma)$ is a discrete series representation of $\tilde U(q+1, p)$ such that its Harish-Chandra parameter does not have any entries in $(-k, 0)$. Then $\sigma= D((\mu_{+}^{+}, \mu_{-}^{-}), (\mu^{-}_{+}, \mu_{-}^{+}))$ with the 
$$\ca(\mu_{+}^{+})+ \ca(\mu_{-}^-)=p-1, \qquad \ca(\mu_{+}^{-})+\ca(\mu_{-}^{+})=q, \qquad\ca(\mu_{-}^{+} \cap [-k, 0]) = \ca(\mu_{-}^{-} \cap [-k, 0])=0.$$
In addition either we have
$$ M_{A(0, -k)}(\sigma)=  D((\mu_{-}^+, -k, \mu_{+}^-), (\mu_{-}^-, 0, \mu^+_{+}))^* \otimes {\det}^{-k}  $$ with
 $\ca(\mu_{+}^{+})+ \ca(\mu_{-}^+)=r-1$ and $\ca(\mu_{+}^{-})+\ca(\mu_{-}^{-})=s;$ or
 $$  M_{A(0, -k)}(\sigma)= D((\mu_{-}^+, 0, \mu_{+}^-), (\mu_{-}^-, -k, \mu^+_{+}))^* \otimes {\det}^{-k}  $$ with
 $\ca(\mu_{+}^{+})+ \ca(\mu_{-}^+)=r$ and $ \ca(\mu_{+}^{-})+\ca(\mu_{-}^{-})=s-1$.
 Conversely, if one of these two situations occurs, then $\sigma \boxtimes M_{A(0, -k)}(\sigma) \in A(0,-k)|_{\tilde U(p-1,q) \times \tilde U(q+1, p)}$.

\end{cor}
Now let $\chi^+= (\mu_{+}^{+}, \mu_{-}^{-})$ and $\chi^-=(\mu^{-}_{+}, \mu_{-}^{+})$. Then
$\sigma=D(\chi^+, \chi^-)$ and 
$$\ca(\chi^+ \cap [-k,0])=\ca(\chi^- \cap [-k, 0]).$$
We have either
\begin{equation}\label{ind1}
M_{A(0,-k)}(\sigma)=D((\chi^{-}, -k),(\chi^{+}, 0))^* \otimes {\det}^{-k};
\end{equation}
or 
\begin{equation}\label{ind2}
M_{A(0,-k)}(\sigma)=D((\chi^{-}, 0), (\chi^{+}, -k))^* \otimes {\det}^{-k}.
\end{equation}
Clearly, up to a ${\det}^m$-character, every discrete series representation of $\tilde U(q+1,p)$ appears as the multiplicity space $M_{A(0,-k)}(\sigma)$ for some $\sigma$.
This will provide us the basis for proving the GGP conjecture inductively and for constructing discrete series representation inductively.

\section{The Law of Reciprocity and GGP Interlacing Relation}
Recall from introduction that if $\pi$ is a unitary representation of $H_1 \times H_2$ and $\sigma_i \in \hat{H_i}$, we have a canonical isometry:
$$\Hom_{H_1}(\sigma_1, M_{\pi}(\sigma_2)) \cong \Hom_{H_2}(\sigma_2, M_{\pi}(\sigma_1)).$$
\\
Suppose that we have the following situation:
$$\begin{array}[c]{ccc}
H_1 &{\subseteq}& G_1\\
\updownarrow\scriptstyle{}&&\updownarrow\scriptstyle{}\\
G_2 &{\supseteq} & H_2 .
\end{array}$$
\begin{enumerate}
\item $(G_1, H_2)$ is a commuting pair of subgroups of $G$;
\item $(G_2, H_1)$ is a commuting pair of subgroups of $G$;
\item $H_1$ is a subgroup of $G_1$ and $H_2$ is a subgroup of $G_2$.
\end{enumerate}
Loosely, this situation is similar to the see-saw dual pairs. \\
\\
Let $\pi$ be a unitary representation of $G$. Suppose that $G_i$ and $H_i$  are CCR. Then $M_{\pi}(\sigma_1)$ is a unitary representation of $G_2$ and $M_{\pi}(\sigma_2)$ is a unitary representation of $G_1$. We have $\sigma_1 \in M_{\pi}(\sigma_2)|_{H_1}$ if and only if $\sigma_2 \in M_{\pi}(\sigma_1)|_{H_2}$. In addition,
$$m( M_{\pi}(\sigma_2)|_{H_1}, \sigma_1)=m( M_{\pi}(\sigma_1)|_{H_2}, \sigma_2).$$
\begin{thm}[Reciprocity] Let $\tau_1 \in \pi|_{G_1}$ and $\sigma_1 \in \hat{H_1}$. Suppose that $M_{\pi}(\tau_1)$ is an irreducible unitary representation of $H_2$ and the related isotypic subspaces
$$\mc H_{\pi}(\tau_1)= \mc H_{\pi}(M_{\pi}(\tau_1))= \tau_1 \hat{\otimes} M_{\pi}(\tau_1).$$
Then $\sigma_1 \in \tau_1|_{H_1}$ if and only if $M_{\pi}(\tau_1) \in M_{\pi}(\sigma_1)|_{H_2}$. In addition, 
$$m(\tau_1|_{H_1}, \sigma_1 )= m(M_{\pi}(\sigma_1)|_{H_2}, M_{\pi}(\tau_1)).$$
\end{thm}
Proof: Take $\sigma_2=M_{\pi}(\tau_1)$. Then $M_{\pi}(\sigma_2) \cong \tau_1$. Our assertion follows immediately.

\subsection{Reciprocity for $A(0, -k)$}
Now we specialize in the following setting:
$$G=\tilde U(p+q, p+q); \qquad \pi=A_{r,s;s,r}(0, -k); $$
$$G_1=\tilde U(p,q); \qquad H_2=\tilde U(q, p); \qquad G_2=\tilde U(q+1, p); \qquad H_1 = \tilde U(p-1, q).$$
We have seen that all conditions in the reciprocity theorem are satisfied. We obtain
\begin{cor}\label{re}[Reciprocity for $A(0,-k)$]  Suppose that $k \geq 1$. Let $\tau_1 \in A(0, -k)|_{\tilde U(p,q)}$ and let $\sigma_1 \in A(0,-k)|_{\tilde U(p-1, q)}$. We have $\sigma_1 \in \tau_1|_{\tilde U(p-1,q)}$ if and only if 
$$M_{A(0,-k)}(\tau_1) \in M_{A(0,-k)}(\sigma_1)|_{\tilde U(q,p)}.$$
 If $\sigma_2 \notin  A(0,-k)|_{\tilde U(q, p)}$, then $\sigma_2 \notin M_{A(0,-k)}(\sigma_1)|_{\tilde U(q,p)}.$   Finally $$m(M_{A(0,-k)}(\sigma_1)|_{\tilde U(q,p)}, M_{A(0,-k)}(\tau_1))= m(\tau_1|_{\tilde U(p-1,q)}, \sigma_1).$$
\end{cor}

\subsection{Harish-Chandra Parameters and GGP Interlacing Relation}
Recall that we define our Harish-Chandra parameter for $\tilde U(p,q)$ to be $(\lambda^+, \lambda^-)$ with $\lambda^+ \in \mathbb R^p$ and $\lambda^- \in \mathbb R^q$. The numerical sequence $\lambda$ then parametrizes the $L$-packets for the discrete series. It is essentially the infinitesimal character. To work with the their conjecture, we adopt the notation of Gan-Gross-Prasad by taking 
$$(\chi, z), \qquad(\chi \in \mathbb R^{p+q}, z \in \mathbb \{\pm 1\}^{p+q})$$
as the Harish-Chandra parameter. Then $\chi^+$ will be the subsequence with $z_i=1$ and $\chi^{-}$ will be the subsequence with $z_i=-1$. We may now write $D(\chi, z)$ for the discrete series corresponding to the Harish-Chandra parameter $(\chi, z)$. We call $z$ the sign sequence. \\
\\
From now on, \uline{we require $\chi$ to be in descending ordering}.
\begin{defn}\label{pattern} Given $(\chi, z)$ and $(\eta, t)$, we mark $z$ by a sequence of $+$ and $-$ and mark  $t$ by a sequence of $\oplus$ and $\ominus$. We line up (the union of) $\chi$ and $\eta$ in a descending ordering. The corresponding sequence of signs extracted from $z$ and $t$ will be called an (interlacing) sign pattern of  $(\chi, z)$ and $(\eta, t)$. We say that $(\chi, z)$ and $(\eta, t)$ satisfy GGP interlacing relation if an interlacing sign pattern of  $(\chi, z)$ and $(\eta, t)$ contains only the following eight consecutive subsequences of size two
$$(\oplus +), (+ \oplus), (- \ominus), (\ominus -), (+-), (-+), (\oplus \ominus), (\ominus \oplus).$$
\end{defn} 
In the case that all entries of $\chi$ and $\eta$ are distinct, there will be only one unique sign pattern. Hence there will be no ambiguity to speak of the (interlacing) sign pattern in this situation. Otherwise, there are ways to line up $\chi$ and $\eta$ in descending ordering. Hence there will be different possible sign patterns. If one of these sign pattern meets the requirement of Def. \ref{pattern}, then we say that $(\chi, z)$ and $(\eta, t)$ satisfy GGP interlacing relation. Our definition is slightly different from the one stated in \cite{gp} \cite{ggp}. We modify it to make the proof of GGP conjecture easier to read.\\
\\
In the case there are no negative signs $\ominus$ and $-$, only $(+ \oplus)$ or $(\oplus +)$ are allowed. Then GGP interlacing relation is the Cauchy interlacing relation:
$$\chi_1 \geq \eta_1 \geq \chi_2 \geq \eta_2 \geq \chi_3 \ldots$$
or
$$  \eta_1 \geq \chi_1 \geq \eta_2 \geq \chi_2 \geq \eta_2 \ldots.$$
If $ \ca(\chi)-\ca(\eta) = \pm 1$, the interlacing relations become  unique. \\
\\
We make the following observation.
\begin{lem}\label{signinduction} Suppose that $(\chi, z)$ and $(\eta, t)$ satisfy the GGP interlacing relation. If their sign pattern has any of the following adjacent pairs
$$(+ \oplus), \qquad (\oplus +), \qquad (- \ominus), \qquad (\ominus -),$$
then by deleting one such pair from the sign pattern and the corresponding entries in $\chi$ and $\eta$, the new pair $(\chi^0, z^0)$ and $(\eta^0, t^0)$ still satisfy the GGP interlacing relation.
\end{lem}
We now state the GGP conjecture that is a reformulation of the original local GGP conjecture (See Ch 17 and 20 of  \cite{ggp}, or  \cite{gp}). In the introduction, we remark that our reformulation is  equivalent what is stated as Conjecture 12.27 in \cite{gp} for unitary groups.\\
\\
{\bf Gan-Gross-Prasad Conjecture}: 
Let $D(\chi,z)$ be a discrete series representation of $U(q,p)$ and $D(\eta, t)$ be a discrete series representation of $U(q-1,p)$. Then $D(\eta, t) \in D(\chi,z)|_{U(q-1,p)}$ if and only if $(\eta, t)$ and $(\chi, z)$ satisfy the GGP interlacing relation and each of these $D(\eta, t)$ has multiplicity one.\\
\\
As we have seen, only discrete series can appear in the square integrable representation $D(\chi, z)|_{U(q-1,p)}$. Hence GGP conjecture gives a complete description of the discrete spectrum of $D(\chi, z)|_{U(q-1,p)}$. We shall remark that the multiplicity one theorem has already been proved in \cite{sz} for all irreducible unitary representations of $U(q,p)$. In this paper, we shall establish the multiplicity one part of the GGP conjecture inductively using the reciprocity theorem. This is different from the proof of Sun and Zhu in nature.\\
\\
The GGP conjecture for $(U(q,p), U(q, p-1))$ can be stated in exactly the same way.  
We observe that if $D(\chi, z)$ is a discrete series of $U(p,q)$, then $D(\chi, -z)$ is a discrete series of $U(q, p)$. If we identify $U(q,p)$ with $U(p,q)$ by identifying the Hermitian form for $U(q,p)$ with the negative of the Hermitian form for $U(p,q)$, then $D(\chi, -z)$ will be equivalent to $D(\chi, z)$. The GGP conjecture for $(U(q,p), U(q, p-1))$ is equivalent to the GGP conjecture for $(U(p,q), U( p-1, q))$.\\
\\
As we have seen, the Harish-Chandra parameters will be either half-integral or  integral for the unitary group $U(q,p)$ depending on the parity of $p+q$. In addition, all entries of Harish-Chandra parameters are distinct. Hence for the GGP conjecture, all the entries of $\chi$ and $\eta$ will be distinct.  There is thus no ambiguities how the two Harish-Chandra parameters should line up in descending ordering. We can then speak of the sign pattern for $(\chi, z)$ and $(\eta, t)$.
\begin{lem}\label{head} If the Harish-Chandra parameters $(\chi, z)$ and $(\eta, t)$ in the GGP conjecture
satisfy the GGP interlacing relation, the sign pattern must start with either with $+$ for $(\chi_1, +1)$ or $\ominus$ for $(\eta_1, -1)$.
\end{lem}
Proof: This Lemma is obviously true for $p=0$ or $q=1$. Applying induction on $p$ and $q$, our assertion is then a consequence of Lemma \ref{signinduction}. $\Box$\\
\\
This lemma is not generally true for all GGP interlacing relations. The key fact we use is that $(\chi, z)$ comes from $U(q,p)$ and $(\eta, t)$ comes from $U(q-1, p)$. This lemma will be crucial in our proof of the GGP conjecture, yet not so obvious from the definition of the GGP interlacing relation.

\section{Proof of Gan-Gross-Prasad Conjecture for $U(p,q)$}

We need the following lemma to prepare for our proof.
\begin{lem}
Let $H$ be a subgroup of $G$. Let $\sigma \in \hat H$ and $\pi \in \hat G$. Then
$\sigma \in \pi|_{H}$ if and only if $\sigma^* \in \pi^*|_{H}$.
\end{lem}
Clearly, $\sigma \in \pi|_{H}$ if and only if
$\sigma^c \in \pi^c |_{H}$. This lemma follow easily from the fact that $\sigma^* \cong \sigma^c$ and $\pi^* \cong \pi^c$. 
\begin{lem}\label{character}
Let $H$ be a subgroup of $G$. Let $\sigma \in \hat H$ and $\pi \in \hat G$. Let $\chi$ be a unitary character of $G$. Then
$\sigma \in \pi|_{H}$ if and only if $\sigma \otimes [\chi|_{H}] \in [\pi \otimes \chi]|_{H}$.
\end{lem}
We shall also need a restatement of Theorem \ref{discrete:equalcase} in terms of $(\chi, z)$. Here we switch $p$ with $q$. 
\begin{prop}\label{discrete:equalcase2} Let $p+q=r+s=n$. Then the discrete spectrum of $A_{r,s;s,r}(0, -k)|_{\tilde U(p,q) \tilde U(q, p)}$ is the direct sum of
$$[D(\eta, -t^{\prime})^* \otimes {\det}^{-k}] \boxtimes D(\eta, t)$$
where
\begin{eqnarray}
& \#\{ \eta_i >0 \mid t_i =1 \}+ \#\{  \eta_i < 0 \mid t_i =-1 \}   =  r, \\
& \#\{\eta_i >0 \mid  t_i =-1 \}+ \#\{  \eta_i < 0  \mid t_i =1 \}  =   s, \\
& \#\{ \eta_i \in (-k, 0) \mid t_i=1 \}  =  \# \{   \eta_i \in (-k, 0) \mid t_i=-1 \},  \\
& t^{\prime}_i  =  \left\{ \begin{array}{cc} 
					t_i & \mbox{ if $\eta_i \notin (-k, 0)$} \\
                  -t_i & \mbox{ if  $\eta_i \in (-k, 0)$}.
                  \end{array}
                  \right.
\end{eqnarray}
All entries in $\eta$ are half integers.            
\end{prop}
Proof: For $\eta_i \notin (-k, 0)$, $\eta_i+k$ remains either positive or negative as $\eta_i$. By the operation of duality between $U(p,q)$ and $U(q, p)$, the corresponding $t_i$ changes sign. For $\eta_i \in (-k, 0)$, $\eta_i+k$ takes the opposite sign of $\eta_i$ as a half integer. By the operation of duality between $U(p,q)$ and $U(q, p)$, the corresponding $t_i$ does not changes sign. $\Box$ \\
\\
Many important theorems in the representation theory of semisimple Lie groups are proved by induction.  Our situation is not any different. We start our induction on $n=p+q$.\\
\\
{\bf Proof of the Gan-Gross-Prasad Conjecture}: 
When $n=2$, we have either $G=U(1,1)$ or $G= U(2)$. In these cases, it is well-known that the GGP conjecture is true.\\
\\
Let us assume that the GGP conjecture is true for any $U(p,q)$ with $p+q=n$. Either we will have $G=U(0, n+1)$, $G=U(n+1, 0)$ or $G=U(q+1, p)$ for $q \geq 0$ and $p \geq 1$.  The irreducible representations of $U(n+1)$ are parametrized by the highest weights.  The branching law has been known to obey the Cauchy interlacing relation and the multiplicity one theorem. In other words, $\sigma_{\mu} \in \pi_{\lambda}|_{U(n)}$ if and only if $\mu$ and $\lambda$ satisfy the Cauchy interlacing relation. In this case, the interlacing relation allows equalities. For compact groups, the Harish-Chandra parameter is the infinitesimal character. It can be obtained from the highest weight by adding the half sum of positive roots of the group. It is not hard to see that $\sigma \in \pi|_{U(n)}$ if and only if the Harish-Chandra parameters of $\pi$ and $\sigma$ satisfy the (strict) Cauchy interlacing relation. As we remarked earlier, when $G$ is compact, the GGP interlacing relation is the Cauchy interlacing relation. Hence the GGP conjecture is true for $U(0, n+1)$ and $U(n+1, 0)$. \\
\\
Now consider $U(q+1, p)$. Here $q \geq 0$ and $p \geq 1$.  
We would like to prove the GGP conjecture for $U(q+1, p)$. Let $(\chi, z)$ be the Harish-Chandra parameter of a discrete series representation of $U(q+1, p)$. Obviously $\#(z_i=1)=q+1$ and $\#(z_i=-1)=p$. Both of these numbers are greater or equal to 1. Consider the sign pattern of $(\chi,z)$.  \,\, \underline{Let the first sign change happen between
$\chi_l$ and $\chi_{l+1}$}.  In other words $l$ is the smallest number such that $z_l z_{l+1}=-1$.
Let $k=\chi_{l}-\chi_{l+1}$. Then $k \geq 1$ is an integer. Consider now $D(\chi, z) \otimes {\det}^{-\chi_{l}}$. In order to do this, we may need to go to the covering of $U(q+1, p)$.
\begin{enumerate}
\item If $p+q+1$ is even, then $p+q$ is odd. Then Harish-Chandra parameters for $U(q+1, p)$ have half integral entries. In this case, we will go to the covering $U(q+1, p)^{\sqrt{\det^{p+q}}}$.
\item If $p+q+1$ is odd, then $p+q$ is even. The Harish-Chandra parameters for $U(q+1, p)$ have integral entries. In this case, we may still go to the covering $U(q+1, p)^{\sqrt{\det^{p+q}}}$ which splits. 
\end{enumerate}
In any case, this setup is  compatible with the double covering involved in Howe's correspondence. \\
\\
Due to Lemma \ref{character}, without loss of generalities, we assume that $\pi=D(\chi, z)$ is a discrete series representation of  $U(q+1, p)^{\sqrt{\det^{p+q}}}$ with
$$\chi_1 > \chi_2 > \ldots > \chi_l (=0) > \chi_{l+1}(=-k) > \chi_{l+2} > \ldots > \chi_{p+q+1}.$$
Let $\pi^{\prime}=D(\chi^{\prime}, z^{\prime})$ where $\chi^{\prime}$ and $z^{\prime}$ are obtained by deleting both the $l$-th entries and the $l+1$-th entries of $\chi$ and $z$. Then $\pi^{\prime}$ is a discrete series representation of $U(q, p-1)^{\sqrt{\det^{p+q}}}$. Now we choose our $(r,s)$ in the following way:
\begin{enumerate}
\item Case (a): If $z_{l}=1$ and $z_{l+1}=-1$, let 
\begin{equation}\label{rs1}
\begin{split}
s= & \#\{ \chi_i^{\prime} < 0 \mid z_i^{\prime}=1 \}, \\
r=  \#\{ \chi_i^{\prime} > 0 \mid z_i^{\prime}=1 \} +  & \# \{ \chi_i^{\prime} < 0 \mid z_i^{\prime}=-1 \}+1=  l+\# \{ \chi_i^{\prime} < 0 \mid z_i^{\prime}=-1 \}.
\end{split}
\end{equation}
Then we have
\begin{eqnarray}
  \# \{ \chi_i^{\prime} < 0 \mid z_i^{\prime}=-1 \}& = & r-l = p-1, \\
  \#\{ \chi_i^{\prime} < 0 \mid z_i^{\prime}=1 \}  & = & s, \\
  \#\{ \chi_i^{\prime} >0 \mid z_i^{\prime} =1 \} & = & q-s=l-1 .
\end{eqnarray}
\item Case (b): If $z_l=-1$ and $z_{l+1}=1$, let
\begin{equation}\label{rs2}
\begin{split}
 r= & \# \{ \chi_i^{\prime} < 0 \mid z_i^{\prime}=-1 \}, \\
 s=  \#\{ \chi_i^{\prime} < 0 \mid z_i^{\prime}=1 \}+1+ &\#\{ \chi_i^{\prime} > 0 \mid z_i^{\prime}=-1 \}=\#\{ \chi_i^{\prime} < 0 \mid z_i^{\prime}=1 \}+l.
 \end{split}
\end{equation}
Then we have
\begin{eqnarray}
\# \{ \chi_i^{\prime} < 0 \mid z_i^{\prime}=-1 \}& =& r,\\
\#\{ \chi_i^{\prime} < 0 \mid z_i^{\prime}=1 \}& = & s-l=q,\\
\# \{ \chi_i^{\prime} > 0 \mid z_i^{\prime}=-1 \} & = & l-1=p-1-r.
\end{eqnarray}
\end{enumerate}
In both cases $p+q=r+s$.
 Denote $A_{r,s; s, r}(0, -k)$ by $A(0, -k)$. By Cor \ref{induction1}, we have
$$[D(\chi^{\prime}, -z^{\prime})^* \otimes {\det}^{-k}] \boxtimes D(\chi, z) \in A(0, -k)|_{\tilde U(p-1, q) \tilde U(q+1, p)}.$$
We apply our Cor. \ref{re}, Equations \ref{ind1} and \ref{ind2}, Prop. \ref{discrete:equalcase2}. Then $D(\eta, t) \in D(\chi, z)|_{\tilde U(q, p)}$ { if and only if} the following are true
\begin{enumerate}
\item $D(\eta, t) \in A(0, -k)|_{\tilde U(q, p)}$. This is equivalent to equalities in Prop \ref{discrete:equalcase2}. Let $d=\#\{  \eta_i >0 \mid t_i =-1 \}$. Then 
$\#\{  \eta_i < 0 \mid  t_i =-1 \}=p-d$ and 
$$\#\{ \eta_i < 0 \mid  t_i =1 \}=s-d, \qquad \#\{  \eta_i >0 \mid t_i =1 \}=q-s+d;$$
$$\#\{  \eta_i \in (-k, 0) \mid t_i=1 \}= \# \{  \eta_i \in (-k, 0) \mid  t_i=-1 \}.$$
All the counting numbers here must be nonnegative and $\eta_i$'s are half integers.
\item $D(\chi^{\prime}, - z^{\prime})^* \otimes {\det}^{-k} \in D(\eta, -t^{\prime})^* \otimes {\det}^{-k}|_{\tilde U(p-1,q)}$ which is equivalent to 
$$D(\chi^{\prime}, -z^{\prime}) \in D(\eta, -t^{\prime})|_{\tilde U(p-1, q)}.$$
Here $ t^{\prime}$ is defined as in Prop. \ref{discrete:equalcase2}.
\end{enumerate}
\vspace{0.2 in}
Suppose that $D(\eta, t) \in D(\chi, z)|_{\tilde U(q, p)}$. Then $D(\chi^{\prime}, -z^{\prime}) \in D(\eta, -t^{\prime})|_{\tilde U(p-1, q)}$. By the induction hypothesis, $(\eta, -t^{\prime})$ and $(\chi^{\prime}, -z^{\prime})$ must satisfy the GGP interlacing relation. Again we use $\oplus$, $\ominus$ for the sign patterns of $\eta$,  and $+$, $-$ for the sign patterns of $\chi$ and $\chi^{\prime}$.  By Lemma \ref{head}, the sign pattern of $(\eta, -t^{\prime})$ and $(\chi^{\prime}, -z^{\prime})$ must start with $\oplus$ for $(\eta_1, +1)$ or $-$ for $(\chi^{\prime}_1, -1)$. Here $-t^{\prime}$ has one more $+1$ than $-z^{\prime}$. It follows that $(\eta, t^{\prime})$ and $(\chi^{\prime}, z^{\prime})$ also satisfy the GGP interlacing relation and the sign pattern must start with either $\ominus$ for $(\eta_1, -1)$ or 
$+$ for $(\chi^{\prime}_1, +1)$.\\
\\
Suppose that $l>1$. Then $\chi^{\prime}_i=\chi_i ( i \leq l-1)$. Consider the sign pattern of $(\eta, t^{\prime})$ and $(\chi^{\prime}, z^{\prime})$ before $\chi_{l}=0$. Here $\chi_{l}$ is not in the sequence. 
\begin{enumerate}
\item Case (a): If $z_i=z_i^{\prime}=1 (i < l)$, there are $l-1$ $(+)$'s, $q-s+d=l-1+d$ $(\oplus)$'s and $d$ $(\ominus)$'s, before $\chi_{l}=0$. Since we start with either $\ominus$ or $+$, by the GGP interlacing relation, the ending sign before $\chi_{l}=0$ must be $\oplus$. So the sign sequence between $0$ and $-k$  will be
$$\hat{0} \ominus \oplus \ominus \oplus \ldots \ominus \oplus \hat{-k}$$
extracted from $(\eta, t^{\prime})$ alone. This sequence could be empty. The first sign after $\chi_{l+1}=-k$ should be either $+$ or $\ominus$. \\
\\
In any case, go back to the sign pattern of the $(\chi, z)$ and $(\eta, t)$ sequence. The last sign before $\chi_{l}=0$ will again be $\oplus$. After this, the sign pattern will go as follows
$$+ \oplus \ominus \oplus \ominus \ldots \oplus \ominus -.$$
Here $+$ comes from $\chi_l=0$, $-$ comes from $\chi_{l+1}=-k$, what are between come from  flipping the signs from $t^{\prime}$ to $t$. The sign pattern after $\chi_{l+1}=-k$ comes from the sign pattern of  $(\eta, t^{\prime})$ and $(\chi^{\prime}, z^{\prime})$ after $-k$. It will start with either $+$ or $\ominus$. Now we have seen that $(\chi, z)$ and $(\eta, t)$ satisfy the GGP interlacing relation.

\item If $z_i=z_i^{\prime}=-1 (i < l)$, there are $l-1$ $(-)$'s, $q-s+d=-l+d$ $(\oplus)$'s and $d$ $(\ominus)$'s, before $x_{l}=0$.  Here obviously $d \geq l$. Since we must start with  $\ominus$ (no $+$ here) , by the GGP interlacing relation, the ending sign before $\chi_{l}=0$ must be $\ominus$. Hence the sign sequence between $0$ and $-k$  will be
$$ \hat{0} \oplus \ominus \oplus \ominus \ldots  \oplus \ominus \hat{-k} $$
extracted from $(\eta, t^{\prime})$ alone. This sequence could be empty. The first sign after $\chi_{l+1}=-k$ should be either $-$ or $\oplus$. \\
\\
In any case, go back to the sign pattern of the $(\chi, z)$ and $(\eta, t)$ sequence. The last sign before $\chi_{l}=0$ will again be $\ominus$. After this, the sign pattern will go as follows
$$-  \ominus \oplus \ominus \oplus \ldots \ominus \oplus  +.$$
Here $-$ comes from $\chi_l=0$, $+$ comes from $\chi_{l+1}=-k$, what are between come from the flipping the signs from $t^{\prime}$ to $t$. The signs  after $\chi_{l+1}=-k$ come from $(\eta, t^{\prime})$ and $(\chi^{\prime}, z^{\prime})$ after $-k$. It will start with either $-$ or $\oplus$. Now we have seen that $(\chi, z)$ and $(\eta, t)$ satisfy the GGP interlacing relation.
\end{enumerate}
Suppose that $l=1$. Then  $\chi_i^{\prime} < -k$ for all $i$. We again consider the sign pattern of $(\eta, t^{\prime})$ and $(\chi^{\prime}, z^{\prime})$. 
\begin{enumerate}
\item Case (a): $z_1=1$ and $z_2=-1$.  Then by Equation \ref{rs1}
$$s=\#\{ \chi_i^{\prime} < 0 \mid z_i^{\prime}=1 \}=q, \qquad
r=\# \{ \chi_i^{\prime} < 0 \mid z_i^{\prime}=-1 \}+1=p-1+1=p.$$ 
As for $(\eta, t^{\prime})$, we have
$$d=\#\{  \eta_i >0 \mid t_i =-1 \}, \qquad \#\{  \eta_i >0 \mid t_i =1 \}=q-s+d=d.$$
Hence there are $d$ $(\oplus)$'s, $d$ $(\ominus)$'s and no $(+)$'s, no $(-)$'s before $\chi_1=0$. Then sign pattern of $(\chi^{\prime}, z^{\prime})$ and $(\eta, t^{\prime})$ must start with $\ominus$ and goes as follows
$$\overbrace{\ominus \oplus \ominus \oplus \ldots \ominus \oplus}^{2d} \hat{0} \overbrace{\ominus \oplus \ominus \oplus \ldots \ominus \oplus}^{\#\{0 > \eta_i > -k\}} \hat{-k} $$ 
then 
$\ominus$ or $+$ immediately after that. In any case, if we insert $+$ into the position of $\hat{0}$, insert
$-$ into the position of $\hat{-k}$ and switch the signs between from $t^{\prime}$ to $t$, then we see that the sign pattern of $(\chi, z)$ and $(\eta, t)$ must be
$$\overbrace{\ominus \oplus \ominus \oplus \ldots \ominus \oplus}^{2d} + \overbrace{ \oplus \ominus \oplus \ominus \ldots  \oplus \ominus}^{\#\{0 > \eta_i > -k\}} -, $$ 
followed by either $\ominus$ or $+$ and the rest of the sign pattern of $(\chi^{\prime}, z^{\prime})$ and $(\eta, t^{\prime})$ after $-k$. Hence  $(\chi, z)$ and $(\eta, t)$ satisfy the GGP interlacing relation. Notice that if $d=0$ or $\#\{0 > \eta_i > -k\} =0$, our argument is still valid.
\item Case (b): $z_1=-1$ and $z_2 = 1$.  Then by Equation \ref{rs2} we have
$$r=\#\{ \chi_j^{\prime} < 0 \mid z_i^{\prime}=-1 \}=p-1, \qquad s= \# \{ \chi_i^{\prime} < 0 \mid z_i^{\prime}=1 \}+1=q+1.$$
As for $(\eta, t^{\prime})$, we have
$$d=\#\{ \eta_i >0 \mid t_i =-1 \}, \qquad \#\{\eta_i >0 \mid t_i =1 \}=q-s+d=d-1.$$
We must have $d \geq 1$.  Then sign pattern of $(\chi^{\prime}, z^{\prime})$ and $(\eta, t^{\prime})$ must start with $\ominus$ and goes as follows
$$\ominus \overbrace{\oplus \ominus \oplus \ominus \ldots \oplus  \ominus}^{2d-2} \hat{0}  \overbrace{ \oplus \ominus \oplus \ominus \ldots  \oplus \ominus}^{\#\{0 > \eta_i > -k\}} \hat{-k} ,$$
then $\oplus$ or $-$ immediately after that. If we insert $-$ into the position of $\hat{0}$, insert
$+$ into the position of $\hat{-k}$ and switch the signs between from $t^{\prime}$ to $t$, then we see that the sign pattern of $(\chi, z)$ and $(\eta, t)$  must go as follows
$$\ominus \overbrace{\oplus \ominus \oplus \ominus \ldots \oplus  \ominus}^{2d-2} -  \overbrace{  \ominus \oplus \ominus \oplus \ldots  \ominus \oplus }^{\#\{0 > \eta_i > -k\}} +$$
followed by either $\oplus$ or $-$ and the rest of the sign pattern of $(\chi^{\prime}, z^{\prime})$ and $(\eta, t^{\prime})$ after $-k$. Hence
$(\chi, z)$ and $(\eta, t)$ satisfy the GGP interlacing relation. Notice that if $d=1$ or $\#\{0 > \eta_i > - k\} =0$, our argument is still valid.
\end{enumerate}
Conversely, if $(\chi, z)$ and $(\eta, t)$ satisfy the GGP interlacing relation, we can go backwards and show that $(\chi^{\prime}, z^{\prime})$ and $(\eta, t^{\prime})$ satisfy the GGP interlacing relation and $ D(\eta, t) \in A(0, -k)|_{\tilde U(q, p)}$. Hence $D(\eta, t) \in D(\chi, z)|_{\tilde U(q, p)}$. \\
\\
The multiplicity one part of the conjecture follows from the induction process and the reciprocity theorem. $\Box$\\
\\
We shall make one remark concerning the proof. One can start the proof based on Theorem \ref{branching2}, rather than Cor. \ref{induction1}. But the proof will be a lot messier. Cor \ref{induction1} allows us to make the $\chi$-part of the interlacing relation rather clean. Thus the proof is mostly concerned with sorting out the $\eta$-part. 
\section{Inductive Construction of $A_{\f q}(\lambda)$ and Branching Laws}

We can now extract the following theorem from the proof of GGP. It is simply the combination of Cor \ref{induction1} with Thm. \ref{multi}.
\begin{prop} Let $\pi=D(\chi, z)$ is a discrete series representation of  $U(q+1, p)^{\sqrt{\det^{p+q}}}$ with integral $\chi:$
$$\chi_1 > \chi_2 > \ldots > \chi_l (=0) > \chi_{l+1}(=-k) > \chi_{l+2} > \ldots > \chi_{p+q+1}.$$
Let $\pi^{\prime}=D(\chi^{\prime}, z^{\prime})$ where $\chi^{\prime}$ and $z^{\prime}$ are obtained by deleting  both the $l$-th entries and the $l+1$-th entries of $\chi$ and $z$. Then
$$D(\chi,z) \cong A_{r,s;s,r}(0, -k)^{\infty} \otimes_{\tilde U(p-1,q)} [ D(\chi^{\prime}, - z^{\prime}) \otimes {\det}^k],$$
where $r,s$ are given by Equation \ref{rs1} and Equation \ref{rs2}. 
\end{prop}
By twisting $D(\chi, z)$,$A_{r,s; s,r}(0, -k)$ and $ D(\chi^{\prime}, - z^{\prime})$ by certain $\det$ characters simultaneously, we obtain
\begin{thm}
All discrete series representations of $U(p,q)$ can be constructed inductively by invariant tensor products with various $A_{r,s;s, r}(k_1, k_2)$ from an irreducible unitary representation of a compact group $U(|p-q|)$.
\end{thm}
It is perhaps a good place to start discussing inductive construction of discrete series for other groups. The important $A_{\f q}(\lambda)$ are the \lq\lq middle dimensional \rq\rq ones. For the group $Mp_{2n}(\mathbb R)$, there are $A_{\f q}(\lambda)$ with Levi group $\tilde U(p,q): \, (p+q=n)$ and $\lambda$ being a character of $\tilde U(p,q)$. We may write them as $A_{p,q}(k)$. Depending on whether $k$ is an integer or half integer and whether $n$ is even or odd, one obtain a unitary representation of either $Sp_{2n}(\mathbb R)$ or $Mp_{2n}(\mathbb R)$.  For the group $O^*(2n)$, again there are $A_{\f q}(\lambda)$ with Levi group $ U(p,q): \, (p+q=n)$ and $\lambda$ being a character of $ U(p,q)$. We may write them as $A_{p,q}(k)$.
Then we can state the following conjecture.
\begin{conj}
All discrete series representation of $Mp_{2n}(\mathbb R)$ can be constructed inductively by invariant tensor product 
$$\mc H_{\sigma} \otimes_{Mp_{2n-2}(\mathbb R)} A_{p,q}(k)^{\infty}$$
for some $k$ and $p,q$ with $p+q=2n-1$, where $\sigma$ is a discrete series representation of $Mp_{2n-2}(\mathbb R)$ .
All discrete series representation of $O^*(2n)$ can be constructed inductively by invariant tensor product 
$$\mc H_{\sigma} \otimes_{O^*(2n-2)} A_{p,q}(k)^{\infty}$$
for some $k$ and $p,q$ with $p+q=2n-1$ where $\sigma$ is a discrete series representation of $O^*(2n-2)$ .
\end{conj}
 For the group $Sp(p,q)$ and $O(2p, 2q)$,  we may consider $A_{\f q}(\lambda)$ with the Levi factor $U(p,q)$. Inductive construction of discrete series in these situations are much more delicate and require a lot more work.
 Nevertheless, all these middle dimensional $A_{\f q}(\lambda)$'s are expected to have similar properties as in the $U(p,q)$-case. More specifically, let $(G(V_1), G(V_2))$ be diagonally embedded in $G(V_1 \oplus V_2)$. Consider $A_{p,q}(k)$ of $G(V_1 \oplus V_2)$. Then the discrete spectrum of
 $A_{p,q}(k)|_{G(V_1) G(V_2)}$ is expected to have multiplicity one and satisfy the hypothesis of the reciprocity theorem. \\
 \\
 {\bf Problem}: Determine the discrete spectrum of $A_{p,q}(k)|_{G(V_1) \times  G(V_2)}$.\\
 \\
 By solving this problem and modifying our approach, one may be able to prove the local GGP conjecture for some other real groups inductively. Currently, the main obstacle towards generalizing our approach to other groups is that one can only obtain the unipotent $A_{p,q}(k)$ through theta lifting of  one dimensional representations. For other $A_{p,q}(k)$, a geometric or analytic realization is required before any branching law can be established. We shall mention the recent work by Savin that treats the branching law of $A_{\f q}(\lambda)|_{Sp_2(\mathbb R) \times SP_{2}(\mathbb R)}$ for $Sp_4(\mathbb R)$. See the appendix in \cite{gg}.\\
\\
 Finally, the inductive construction is expected to be potentially more general.
\begin{conj} All $A_{\f q}(\lambda)$'s for $U(p,q)$, $Mp_{2n}(\mb R)$ and $O^*(2n)$ can be constructed inductively from an irreducible finite dimensional unitary representation by invariant tensor products with various $A_{p,q}(\lambda)$. 
\end{conj}
The reader should perhaps consult our other papers \cite{hen} \cite{hearthur} to see how invariant tensor product can be used to construct unipotent representations.



\begin{thebibliography}{99}





\bibitem[BW]{bw}  A. Borel, N. Wallach {\it Continuous cohomology, discrete subgroups, and representations of reductive groups,}, Annals of Math. Studies, Princeton University Press, Princeton, (1980).
\bibitem[BP]{bp} R. Beuzart-Plessis \lq\lq A local trace formula for the Gan-Gross-Prasad
conjecture for unitary groups:  the archimedean case, \rq\rq arXiv:1506.01452v2.

\bibitem[Di]{di} J. Dixmier {\it $C^*$-algebra}, North-Holland Publishing Company (1977).

\bibitem[FD]{fd}  J. M. G. Fell, R. S. Doran, {\it Representations of $*$-algebras,Locally Compact Groups, and Banach *-algebraic Bundles} Academic Press (1988).

\bibitem[GGP]{ggp} W.-T. Gan, B.  Gross, D. Prasad \lq\lq Symplectic local root numbers, central critical L values, and restriction problems in the representation theory of classical groups. Sur les conjectures de Gross et Prasad. I.\rq\rq  Asterisque No. 346 (2012), 1-109.
\bibitem[GG]{gg} W.-T. Gan, N. Gerivich \lq\lq Restrictions of Saito-Kurokawa representations,\rq\rq
with an appendix by Gordan Savin. Contemp. Math., 488, Automorphic forms and L-functions I. Global aspects, 95-124, Amer. Math. Soc., Providence, RI, 2009.
\bibitem[GP]{gp} B. Gross, D. Prasad \lq\lq On the decomposition of a representation of $SO_n$ when restricted to $SO_{n-1}$ \rq\rq, {\it Canad. J. Math. } No. 5 Vol. 44 (1992), 974-1002.
\bibitem[GW]{gw} B. Gross and N. Wallach \lq\lq Restriction of small discrete series representations to symmetric subgroups.\rq\rq The mathematical legacy of Harish-Chandra (Baltimore, MD, 1998), 255?272, Proc. Sympos. Pure Math., 68, Amer. Math. Soc., Providence, RI, (2000), 255-272.
\bibitem[MH]{mh} M. Harris \lq\lq Testing rationality of coherent cohomology of Shimura varieties.\rq\rq  Automorphic forms and related geometry: assessing the legacy of I. I. Piatetski-Shapiro,
Contemp. Math., 614, Amer. Math. Soc., Providence, RI, 2014,  81-95.
\bibitem[He00]{theta} H. He \lq\lq  Theta 
Correspondence I-Semistable Range: Construction and Irreducibility 
\rq\rq,  {\it Communications in Contemporary Mathematics } (Vol. 2), (2000), 
255-283.
\bibitem[Heq]{basic} H. He, \lq\lq Composition of Theta Correspondences \rq\rq, {\it Adv. in Math. 190}, (2005), 225-263.
\bibitem[Heu]{unit} H. He, \lq\lq  Unitary Representations and Theta Correspondence for Type I Classical Groups,\rq\rq Journal of Functional Analysis, Vol 199, Issue 1, (2003), 92-121. 
\bibitem[Hen]{hen} H. He, \lq\lq Unipotent Representations and Quantum Induction, \rq\rq preprint (2005).
\bibitem[HA]{hearthur} H. He, \lq\lq Certain Unitary Langlands-Vogan Parameter for Special Orthogonal Groups\rq\rq, preprint (2011).
\bibitem[Howe]{howe} R. Howe, \lq \lq Transcending Classical Invariant 
Theory\rq \rq
{\it J. of Amer. Math. Soc.} Vol. 2, (1989), 535-552.
\bibitem[Ko]{ko} T. Kobayashi \lq\lq Discrete decomposability of the restriction of
$A_{\f q}(\lambda)$
with respect to reductive
subgroups and its applications \rq\rq, Inventiones Mathematicae
117 No. 2,
(1994),  181-205.
\bibitem[Ko2]{ko2} T. Koyayashi \lq\lq Branching problems of Zuckerman derived functor modules \rq\rq {\it Representation Theory and Mathematical Physics}, American Mathematical Society, Providence (2011) 23-40.
\bibitem[KV]{kv} A. Knapp, D. Vogan {\it Cohomological Induction and Unitary Representations}, Princeton University Press, Princeton, NJ 1995.
\bibitem[Ku]{ku} S. Kudla, \lq\lq On the Local Theta Correspondence, \rq\rq {\it Invent. Math.}, V. 83, (1986), 229-255.
\bibitem[Li]{li} J-S. Li, \lq \lq Singular Unitary Representation of 
Classical Groups\rq \rq
{\it Inventiones Mathematicae} V. 97 (1989) 237-255.
\bibitem[Li1]{li1} J.-S. Li, \lq\lq Theta Lifting for Unitary Representations with Nonzero Cohomology \rq\rq, {\it Duke Mathematical Journal}, No. 3 V. 61, (1990), 913-937.
\bibitem[MW]{mw} Moeglin, C.; Waldspurger, J.-L. Sur les conjectures de Gross et Prasad. II. (French) [On the conjectures of Gross and Prasad. II] Asterisque No. 347 (2012).
\bibitem[P1]{p1} A. Paul, \lq\lq Howe Correspondence for Real Unitary Groups,\rq\rq, {\it Journal of Functional Analysis}, No 158, (1998), 384-431.
\bibitem[P2]{p2} A. Paul, \lq\lq Howe Correspondence for Real Unitary Groups II,\rq\rq, {\it Proc. of  Amer. Math. Soc.}, No 10, Vol 128, (2000), 3129-3136.
\bibitem[PT]{pt} A. Paul, P. Trapa, \lq\lq One dimensional Representation of $U(p,q)$ and Howe Correspondence,\rq\rq {\it J. of Func. Analysis}, V. 195, 2002, (129-166).
\bibitem[Sch]{sch} W. Schmid, \lq\lq Discrete Series,\rq\rq {\it Representation Theory and Automorphic Forms},  noted by V. Bolton (83-113), Proc. Sympos. Pure Math., 61, Amer. Math. Soc., Providence, RI, 1997.
\bibitem[SZ]{sz} B. Sun and C. Zhu, \lq\lq Multiplicity One theorems: the Archimedean case.\rq\rq
{\it Ann. of Math.} No 1 Vol 175, 2012 (23-44).
\bibitem[V]{v} D. Vogan \lq\lq Unitarizability of certain series of representations, \rq\rq {\it Annals of Math.} V 120, 1984, (141-187).
\bibitem[VZ]{vz} D. Vogan and G. Zuckerman, \lq\lq Unitary Representations with nonzero cohomology,\rq\rq, {\it Compositio Math.}, V. 53, 1994, (51-90).
\bibitem[WA]{wallach} N. Wallach, {\it Real Reudctive Groups I II}, Academic Press, San Diego CA, 1992.
\bibitem[W]{w}  J.-L. Waldspurger,  Une formule intégrale reli\'{e}e \'{a} la conjecture locale de Gross-Prasad.  Compos. Math. 146 (2010), no. 5, 1180-1290.
\bibitem[Zhang]{zh} W. Zhang, \lq\lq Fourier transform and the global Gan-Gross-Prasad
conjecture for unitary groups,\rq\rq {\it Annals of Math}, Vol 180, (2014) 971-1049.





\end{thebibliography}
\end{document}